# MODULAR CURVES OF GENUS 2

ENRIQUE GONZÁLEZ-JIMÉNEZ AND JOSEP GONZÁLEZ

ABSTRACT. We prove that there are exactly 149 genus two curves $C$ defined over $\mathbb{Q}$ such that there exists a nonconstant morphism $\pi: X_1(N) \to C$ defined over $\mathbb{Q}$ and the jacobian of $C$ is $\mathbb{Q}$-isogenous to the abelian variety $A_f$ attached by Shimura to a newform $f \in S_2(\Gamma_1(N))$. We determine the corresponding newforms and present equations for all these curves.

## 1. INTRODUCTION

Let $C$ be a nonsingular projective curve defined over $\mathbb{Q}$. If the genus of $C$ is 1, the modularity[1] of $C$ as a curve is equivalent to the modularity of $C$ as an abelian variety. That is, for a curve of genus 1 is equivalent the existence of a nonconstant morphism from the modular curve $X_1(N)$ to $C$ defined over $\mathbb{Q}$ to the existence of a morphism from the jacobian of $X_1(N)$, $J_1(N)$, onto the jacobian of $C$, $J(C)$, defined over $\mathbb{Q}$. When the genus is greater than 1, this equivalence does not hold.

This paper is devoted to determine a certain class of hyperelliptic modular curves of genus 2 defined over $\mathbb{Q}$. Several hyperelliptic modular curves are known and their equations have been computed. In [Ogg74] the values of $N$ such that $X_0(N)$ is hyperelliptic are determined and equations for each one of these curves are given in [Gon91]. In [Mes81] it is proved that $N = 13, 16$ and $18$ are the only values such that $X_1(N)$ is hyperelliptic. The values of $N$ such that $X_0^*(N) := X_0(N)/B(N)$ ($B(N)$ denotes the group of the Atkin-Lehner involutions) is hyperelliptic are determined in [HH96] and [Has97], and equations for these curves are presented. In [FH99] the hyperelliptic curves obtained as quotient of $X_0(N)$ by proper subgroups of $B(N)$ are given. More literature about this topic can be also found in [IM91], [Mur92] and [Shi95].

In section 2, we introduce the notion of primitive modular curve of level $N$ over a number field $K$ as a curve $C$ which can be parametrized from $X_1(N)$ over $K$ and such that for any proper divisor $d|N$ the jacobian of $X_1(d)$, $J_1(d)$, does not contain the jacobian of $C$, $J(C)$, as $K$-factor, up to $K$-isogeny. We prove that there are a finite number of such curves $C$ over $\mathbb{Q}$ of genus two with the additional requirement that their jacobians are $\mathbb{Q}$-simple, i.e., for every such curve $C$ there exists a newform $f \in S_2(\Gamma_1(N))$ such that $J(C)$ is $\mathbb{Q}$-isogenous to $A_f$, where $A_f$ is the abelian variety attached by Shimura [Shi71] to $f$.

Date: September 22, 2000.
2000 *Mathematics Subject Classification.* Primary 14G35, 14H45; Secondary 11F11, 11G10.
The first author was supported in part by DGICYT Grant PB96-1154.
The second author was supported in part by DGICYT Grant BFM2000-0794-C02-02.

[1]A recent joint work of C. Breuil, B. Conrad, F. Diamond and R. Taylor [BCDT99] completes Wiles' results to prove the *Conjecture of Shimura-Taniyama-Weil* so that it establishes that every elliptic curve defined over $\mathbb{Q}$ is modular.

The rest of the sections are devoted to finding the corresponding newforms and equations for these curves. In section 3, we establish a program to compute a set of hyperelliptic curves of genus 2 which contains the set of the mentioned curves. In Section 4 we discard the curves obtained in section 3 that are not primitive modular with $\mathbb{Q}$-simple jacobian. In section 5 we determine the corresponding newforms of the remaining curves and in section 6 we check that all these are primitive modular with $\mathbb{Q}$-simple jacobian. As a consequence, we prove that there are only 149 cases, whose equations and corresponding newforms are given in the tables in the appendix. We conclude the paper by giving an example of a genus 2 curve $C$ defined over $\mathbb{Q}$ whose jacobian is a $\mathbb{Q}$-simple quotient of the new part of the jacobian $J_1(36)$ and such that $C$ cannot be parametrized from $X_1(36)$.

**Acknowledgements:** We are grateful to M. Müller, J. Quer and W. A. Stein for their computational help which has allowed us to complete the tables in the appendix.

## 2. Primitive modular curves of genus 2

Let $K$ be a number field. We fix the following notation. As usual, an abelian variety $A$ defined over $K$ will be called $K$-simple if does not contain a proper nontrivial abelian subvariety such that it and its inclusion are defined over $K$. For two abelian varieties $A$ and $B$ defined over $K$, the notation $A \stackrel{K}{\sim} B$ means that $A$ and $B$ are $K$-isogenous. Concerning to the modular curves, we will follow the notation of [Rib77] and we denote by $S_2(N,\varepsilon)$ the complex vector space of cusp forms of weight 2 on $\Gamma_1(N)$ with Nebentypus $\varepsilon$.

Next, we give the following two definitions.

**Definition 2.1.** *We say that an abelian variety $A$ defined over $K$ is modular of level $N$ if there exists a surjective morphism $\nu : J_1(N) \to A$ defined over $K$. In that case, we say that $A$ is primitive of level $N$ if $A$ is not modular of level $d$ for any proper divisor $d|N$.*

**Definition 2.2.** *We say that a nonsingular projective curve $C$ defined over $K$ is modular of level $N$ if there exists a nonconstant morphism $\pi : X_1(N) \to C$ defined over $K$. In that case, we say that $C$ is primitive of level $N$ if its jacobian, $J(C)$, is primitive of level $N$.*

We remind the reader that modularity depends on the field of definition. It is possible to have a curve $C/K$ (resp. abelian variety $A/K$) that is not modular, and that by enlarging the field of definition to an extension $L/K$, we obtain a modular curve $C/L$ (resp. abelian variety $A/L$). For instance, the elliptic curve $E : y^2 = x^3 + A\,x + B$, where $A = -135 + 156\sqrt{-3}$ and $B = 82 + 1092\sqrt{-3}$, is a modular curve (or modular abelian variety) over $K = \mathbb{Q}(\sqrt{-3})$ of level $N = 63$. In fact, $E$ is a $\mathbb{Q}$-curve of degree 3 completely defined over $K$, i.e., there is an isogeny of degree 3 from $E$ to its Galois conjugate defined over $K$. The twisted elliptic curve $E_\gamma : y^2 = x^3 + \gamma^2 A\,x + \gamma^3 B$, where $\gamma = 2 + \sqrt{-3}$, is not modular over $K$ for any level $N$ because it is not completely defined defined over $K$ (cf. [GL01]). However, $E_\gamma$ is modular over $K(\sqrt{\gamma})$.

Let $\pi : X_1(N) \to C$ be a nonconstant morphism defined over $K$. Then, the corresponding morphism between the jacobians $\pi_* : J_1(N) \to J(C)$ is defined over $K$. Therefore, if $C/K$ is modular of level $N$, then $J(C)/K$ is also modular of level $N$. The converse is not true in general (see the example of section 7).

We will restrict our attention to the rational case and, from now on, we suppose $K = \mathbb{Q}$. If $J(C)$ is modular of level $N$, then there exist $k$ normalized eigenforms $f_1, \ldots, f_k$ in $S_2(\Gamma_1(N))$ such that

$$J(C) \stackrel{\mathbb{Q}}{\sim} A_{f_1}^{n_1} \times \cdots \times A_{f_k}^{n_k},$$

where $A_{f_i}$ is the abelian variety defined over $\mathbb{Q}$ attached to $f_i$ by Shimura. If $J(C)$ is $\mathbb{Q}$-simple, then $J(C) \stackrel{\mathbb{Q}}{\sim} A_f$ for some eigenform $f \in S_2(\Gamma_1(N))$. Moreover, if $J(C)$ is primitive of level $N$, then $f$ is a newform and we have the following equality of complex vector spaces:

$$\nu^*(H^0(J(C), \Omega^1)) = \langle {}^\sigma f(q)\, dq/q : \sigma \in \mathrm{Gal}(\overline{\mathbb{Q}}/\mathbb{Q}) \rangle,$$

where, as usual, $q = e^{2\pi i z}$. In this case, thanks to the results of Carayol [Car86] the level $N$ is determined by the conductor of $J(C)$ and the newform $f$ is unique, up to Galois conjugation (cf. [Rib77]).

For the sake of simplicity, given a newform $f \in S_2(\Gamma_1(N))$, we will denote by $S_2(A_f)$ the $\mathbb{C}$-vector space generated by the Galois conjugates of $f(q)$ and we will identify $H^0(A_f, \Omega^1)$ with the $\mathbb{C}$-vector space $\{g(q)\, dq/q : g \in S_2(A_f)\}$.

**Proposition 2.1.** *Let $f \in S_2(N, \varepsilon)$ be a newform such that $A_f$ is an abelian surface. Then the following conditions are equivalent:*

(1) *There exists a primitive modular curve $C$ defined over $\mathbb{Q}$ of level $N$ such that $J(C)$ is $\mathbb{Q}$-isogenous to $A_f$.*
(2) *There exists a hyperelliptic curve $C'/\mathbb{C}$ and a nonconstant morphism $\pi' : X_1(N) \to C'$ defined over $\mathbb{C}$ such that $\pi'^*(H^0(C', \Omega^1)) = H^0(A_f, \Omega^1)$.*
(3) *For every linearly independent pair of modular forms $g_1, g_2 \in S_2(A_f)$, there exists $P(X) \in \mathbb{C}[X]$ of degree 5 or 6 without double roots such that the functions on $X_1(N)$ given by*

$$x = \frac{g_1}{g_2}, \quad y = \frac{q\, dx/dq}{g_2}$$

*satisfy the equation $y^2 = P(x)$.*

*Proof.* It is obvious that (1) implies (2). Let us assume (2) and we will prove (3). Let $Q(X)$ be a polynomial in $\mathbb{C}[X]$ of degree 5 or 6, without double roots and such that $V^2 = Q(U)$ is an affine model for $C'$. We denote by $u$ and $v$ respectively the functions on $C'$ corresponding to the affine coordinates $U$ and $V$ resp. Since the differentials $u\, du/v$, $du/v$ form a basis of the $\mathbb{C}$-vector space $H^0(C', \Omega^1)$, the regular differentials

$$\omega_1 = \pi'^*\left(u\frac{du}{v}\right), \qquad \omega_2 = \pi'^*\left(\frac{du}{v}\right)$$

are linearly independent and the functions on $X_1(N)$:

$$x := \pi'^*(u) = \frac{\omega_1}{\omega_2}, \quad y := \pi'^*(v) = \frac{\pi'^*(du)}{\omega_2} = \frac{d\,x}{\omega_2},$$

satisfy $y^2 = Q(x)$. Notice that we have

$$x = \frac{h_1}{h_2} \quad \text{and} \quad y = \frac{q\,dx/dq}{h_2},$$

where $h_1$, $h_2$ are the modular forms $h_1 = \omega_1 \, q/dq$, $h_2 = \omega_2 \, q/dq$. Now, if we take two linearly independent $g_1, g_2 \in S_2(A_f)$, then there exists a matrix $\begin{pmatrix} a & b \\ c & d \end{pmatrix} \in GL_2(\mathbb{C})$ such that
$$\begin{cases} g_1 = a \, h_1 + b \, h_2 \\ g_2 = c \, h_1 + d \, h_2 \,. \end{cases}$$
Now, we take the polynomial
$$P(X) = \frac{(-cX+a)^6}{(a\,d - b\,c)^4} \, Q\left(\frac{dX - b}{-cX + a}\right) \in \mathbb{C}[X] \,.$$
Let us denote by $C''$ the normalization of the projective curve given by the affine model $V^2 = P(U)$ and by $\phi : C' \to C''$ the isomorphism given by:
$$(u, v) \mapsto (U, V) = \left(\frac{au + b}{cu + d}, v \, \frac{ad - bc}{(cu + d)^3}\right) \,.$$
It is easy to check that the functions on $X_1(N)$
$$x := (\phi \circ \pi')^*(U) = \frac{g_1}{g_2}, \quad y := (\phi \circ \pi')^*(V) = \frac{qd(\phi \circ \pi')^*(U)/dq}{g_2}$$
satisfy the equation $y^2 = P(x)$.

To end the proof, we have to prove that (3) implies (1), but this follows easily from the fact that $S_2(A_f)$ contains a basis $\{h_1, h_2\}$ with rational $q$-expansion and, therefore, its corresponding polynomial $P(X)$ lies in $\mathbb{Q}[X]$.  □

**Remark 2.1.** *Given a hyperelliptic curve $C/\mathbb{C}$ of genus $g$ ($g \geq 2$) defined by an equation $y^2 = P(x)$, where $P(x)$ is a complex polynomial without double roots, the group $G := GL_2(\mathbb{C})$ provides isomorphisms between $C$ and all the isomorphic hyperelliptic curves given by similar equations by means of the following map:*
$$\begin{pmatrix} m & n \\ p & q \end{pmatrix} \in GL_2(\mathbb{C}) \quad \rightsquigarrow \quad (x, y) \mapsto \left(\frac{mx + n}{px + q}, \frac{m\,q - n\,p}{(px + q)^{g+1}} y\right) \,.$$
*Only when the genus of $C$ is two does the group $G$ coincides with the group of the automorphisms of the regular differentials on the curve $C$. This fact, which has been used in the previous proof, cannot be applied when the genus of $C$ is $> 2$.*

¿From now on, let $f \in S_2(N, \varepsilon)$ be a newform with $q$-expansion
$$f = \sum_{n \geq 1} a_n \, q^n \,.$$

Notice that the coefficients $a_n$ are algebraic integers and satisfy the recursive relations:

$$a_{mn} = a_m a_n \quad \text{if } (m, n) = 1 \,, \quad a_{p^k} = a_{p^{k-1}} a_p - p \, \varepsilon(p) a_{p^{k-2}} \quad \text{for all primes } p \,.$$

Furthermore, we have $|a_p| \leq 2\sqrt{p}$ and if $p|N$ then:
$$|a_p| = \begin{cases} 0 & \text{if } p^2|N \text{ and } \mathfrak{f}_\varepsilon \mid N/p, \\ \sqrt{p} & \text{if } \mathfrak{f}_\varepsilon \nmid N/p, \\ 1 & \text{if } p^2 \nmid N \text{ and } \mathfrak{f}_\varepsilon \mid N/p, \end{cases}$$

where $\mathfrak{f}_\varepsilon$ is the conductor of $\varepsilon$ (cf. [DS74]). In particular, if $\varepsilon = 1$, then $a_p$ only takes the values $0, \pm 1$.

In the sequel, we assume that $\dim A_f = 2$. Let $d$ be the square-free integer such that $\mathbb{Q}(\sqrt{d})$ is the field $K_f = \mathbb{Q}(\{a_n\}_{n\geq 1})$ and denote by $\sigma$ the nontrivial automorphism in $\mathrm{Gal}(K_f/\mathbb{Q})$. Let $\{h_1, h_2\}$ be the basis of the $\mathbb{C}$-vector space $\langle f, {}^\sigma f \rangle$ given by:

$$
(1) \quad \begin{cases} h_1 := \dfrac{f + {}^\sigma f}{2} = q + \displaystyle\sum_{n\geq 2} \dfrac{a_n + {}^\sigma a_n}{2} q^n \\ h_2 := \dfrac{f - {}^\sigma f}{2\sqrt{d}} = \displaystyle\sum_{n\geq 2} c_n\, q^n \,. \end{cases}
$$

We have that $h_1, h_2 \in 1/2\,\mathbb{Z}[[q]]$. Let us define

$$n_0 := \min\{n \in \mathbb{Z} : c_n \neq 0\} = \min\{n \in \mathbb{Z} : a_n \notin \mathbb{Z}\}.$$

**Lemma 2.2.** *If $C$ is a primitive modular curve of level $N$ of genus $2$ such that $J(C) \stackrel{\mathbb{Q}}{\sim} A_f$, then $n_0 = 2$ or $n_0 = 3$.*

*Proof.* We take:

$$
(2) \quad \begin{cases} x = \phantom{aaa} c_{n_0} \dfrac{h_1}{h_2} \phantom{aa} = \dfrac{1}{q^{n_0-1}} + \cdots, \\ y = \dfrac{c_{n_0}}{1 - n_0} \dfrac{q}{h_2} \dfrac{dx}{dq} = \dfrac{1}{q^{2n_0-1}} + \cdots. \end{cases}
$$

Now, $x$ and $y$ satisfy $y^2 = P(x)$, where $P(x)$ is a polynomial of degree 5 or 6. As a result of this relation, we have:

$$\mathrm{ord}_{i\infty} y^2 = \mathrm{ord}_{i\infty} x^{\deg P}.$$

Therefore, if $\deg P = 6$ then $2(2n_0 - 1) = 6(n_0 - 1)$ and $n_0 = 2$, whereas if $\deg P = 5$ then we have $2(2n_0 - 1) = 5(n_0 - 1)$ and therefore $n_0 = 3$. □

**Theorem 2.3.** *There is a finite number of primitive modular curves of genus two over $\mathbb{Q}$, up to $\mathbb{Q}$-isomorphism, such that their jacobians are $\mathbb{Q}$-simple.*

*Proof.* We note that the above statement is equivalent to claiming that the set of normalized newforms $f \in S_2(N, \varepsilon)$ such that $A_f \stackrel{\mathbb{Q}}{\sim} J(C)$ for some primitive modular curve $C$ of level $N$ and genus 2 is finite.

Let $f = \sum_{n\geq 1} a_n\, q^n$, $d$ and $\sigma$ be as before. Let $x, y$ be as in (2) and $P \in \mathbb{Q}[X]$ be such that $y^2 = P(x)$. The polynomial $P$ is determined by the values $a_2, a_3, a_5, a_7$, $a_{11}, a_{13}, \varepsilon(2), \varepsilon(3)$. Therefore, it suffices to prove that there is a finite number of possible values of $d$ and a finite number of possibilities for $a_2, a_3, a_5, a_7, a_{11}, a_{13}, \varepsilon(2)$ and $\varepsilon(3)$.

We shall use the notation $a_n = x_n + y_n \sqrt{d}$ with $x_n, y_n \in 1/2\,\mathbb{Z}$. First, we will see that there exists only a finite number of possible values of $d$. We distinguish two cases according to the value $n_0$ as in the previous lemma is 2 or 3. Indeed, if $n_0 = 2$ then

$$0 < |y_2| \sqrt{|d|} = \dfrac{|a_2 - {}^\sigma a_2|}{2} \leq \dfrac{|a_2| + |{}^\sigma a_2|}{2} \leq 2\sqrt{2}.$$

The minimum value for $|y_2|$ is $1/2$ or $1$ depending on whether $d \equiv 1 \pmod{4}$ or not. Therefore, if $d \equiv 1 \pmod{4}$ then $|d| \leq 31$, otherwise $|d| \leq 7$. In the case $n_0 = 3$ with a similar argument about $a_3$ we obtain that if $d \equiv 1 \pmod 4$ then $|d| \leq 47$, otherwise $|d| \leq 11$. So, we have seen that there is only a finite number of possible values of $d$.

To finalize, we will see that for each $d$ there is only a finite number of possible values for $\varepsilon(2), \varepsilon(3)$ and $a_p$ when $p$ runs over $\{2, 3, 5, 7, 11, 13\}$. Indeed, the only possible nonzero values for $\varepsilon(n)$ are $\pm 1, \pm i, \pm(1 \pm \sqrt{-3})/2$, since these are the only roots of the unity contained in a quadratic field. In particular, $\varepsilon(2)$ and $\varepsilon(3)$ must be one of these values or zero. The values $a_p$ are restricted by the conditions:

$$|x_p| = \frac{|a_p + {}^\sigma a_p|}{2} \leq 2\sqrt{p}, \quad |y_p| = \frac{|a_p - {}^\sigma a_p|}{2\sqrt{|d|}} \leq 2\sqrt{p/|d|}.$$

So, there is a finite number of possible values for $a_2, a_3, a_5, a_7, a_{11}, a_{13}$. □

**Remark 2.2.** *Notice that the previous theorem holds even if we allow that the jacobians of the primitive modular curves $C$ of level $N$ can be $\mathbb{Q}$-isogenous to $A_{f_1} \times A_{f_2}$ with $f_1, f_2$ being newforms of $S_2(\Gamma_1(N))$. Indeed, if this is the case then results similar to those in Proposition 2.1 and Lemma 2.2 can be obtained changing $f$ and ${}^\sigma f$ by $f_1$ and $f_2$. However, this case will not be treated in this paper since we are interested in newforms $f$ such that the corresponding field $K_f$ is a quadratic field.*

**Remark 2.3.** *The bounds given for $d$ in the previous proof can be improved by considering whether $f$ has an extra twist or not. We recall that an extra twist of $f$ is a pair $(\sigma, \chi)$, where $\sigma$ is a nontrivial element in $\mathrm{Gal}\,(\overline{\mathbb{Q}}/\mathbb{Q})$ and $\chi$ is a character of $\mathrm{Gal}\,(\overline{\mathbb{Q}}/\mathbb{Q})$, satisfying ${}^\sigma a_p = \chi(p) a_p$ for all primes $p$ not dividing $N$. It is known that if $\varepsilon = 1$, then $\chi^2 = 1$ and if $\varepsilon$ is not trivial, then $(\sigma, \varepsilon^{-1})$ is an extra twist of $f$, when we take $\sigma$ as the complex conjugation. In particular, if $f$ does not have an extra twist then $\varepsilon = 1$ and $d > 0$.*

*When $f$ has an extra twist $(\sigma, \chi)$ the values $a_p$ and $d$ have to satisfy additional conditions. For instance, if the order of $\chi$ is equal to 2 then for all primes $p \nmid N$:*

$$a_p = \begin{cases} x_p \in \mathbb{Z} & \text{if } \chi(p) = 1, \\ y_p \sqrt{d} \in \mathbb{Z}[\sqrt{d}] & \text{if } \chi(p) = -1, \end{cases}$$

*while if $\mathrm{ord}\,\chi \neq 2$ then the only two possibilities for $d$ are $d = -1$ or $d = -3$ depending on whether the order of $\chi$ is 4 or a multiple of 3 (3 or 6).*

*Finally, we obtain the following table of possibilities for $d$:*

| $f$ | $n_0 = 2$ | $n_0 = 3$ |
|---|---|---|
| With extra twist | $d \in \left\{ \begin{array}{c} -1, \pm 2, \pm 3, \\ \pm 5, \pm 6, \pm 7 \end{array} \right\}$ | $d \in \left\{ \begin{array}{c} -1, \pm 2, \pm 3, \pm 5, \\ \pm 6, \pm 7, \pm 10, \pm 11 \end{array} \right\}$ |
| Without extra twist | $d \in \left\{ \begin{array}{c} 2, 3, 5, 6, 7, \\ 13, 17, 21, 29 \end{array} \right\}$ | $d \in \left\{ \begin{array}{c} 2, 3, 5, 6, 7, 10, 11, 13, \\ 17, 21, 29, 33, 37, 41 \end{array} \right\}$ |

In order to determine the levels and equations for the curves $C/\mathbb{Q}$ as in Theorem 2.3, we proceed as follows:

1. **Collecting.** We have designed a program in Mathematica and PARI/GP to compute the possible primitive modular curves of genus two over $\mathbb{Q}$.

2. **Discarding.** We present sieves to reject the curves obtained in the previous stage whose jacobians are not $\mathbb{Q}$-simple or are not modular. These sieves are based on properties of hyperelliptic curves of genus 2 or of jacobians of modular curves.
3. **Searching.** For every curve $C$ which passes all the sieves we look for the possible newforms $f$ such that $A_f \stackrel{\mathbb{Q}}{\sim} J(C)$.
4. **Checking.** We establish a criterion to check when a genus 2 curve corresponding to a newform $f \in S_2(N, \varepsilon)$ is one of the solutions considered in Theorem 2.3.

## 3. Collecting

To determine the possible equations for the curves $C/\mathbb{Q}$ as in Theorem 2.3, we have designed a program consisting of the following steps:

- We introduce the formal expression $f = q + a_2\, q^2 + \cdots + a_M\, q^M + O(q^{M+1})$, with $a_n = x_n + y_n\sqrt{d}$ and $M \geq 13$. The expression $f$ depends only on the variables $d$ and $x_p, y_p, \varepsilon(p)$ with $p$ running over the set of the primes $\leq M$, and all these variables are subject to the restrictions mentioned in the proof of the Theorem 2.3 and the remark 2.3.
- We write $x$ and $y$ as in (2) and compute the polynomial

$$P(X) = \sum_{i=0}^{6} A_i\, X^i$$

such that $y^2 - P(x) = O(q)$. We ignore $P(x)$ when this polynomial has double roots.
- We compute the polynomial in $q$: $y^2 - P(x) = \sum_{m=1}^{M'} b_m(\{a_n\}_{2 \leq n \leq M})\, q^m$ and solve the system

$$b_m(\{a_n\}_{2 \leq n \leq M}) = 0, \quad m = 1, \ldots, M'$$

in the variables $a_n$, $2 \leq n \leq M$, and where $M' = M - 14$ or $M' = M - 8$ depending on whether the degree of the polynomial $P(x)$ is 5 or 6.
- Finally, we obtain $(P(X), \{a_n\}_{2 \leq n \leq M}, k, d)$, where the curve $C$ is defined by $y^2 = P(x)$, $f = q + a_2\, q^2 + \cdots + a_M\, q^M + O(q^{M+1}) \in S_2(N, \varepsilon)$, $k = \text{ord}\, \varepsilon$ and $K_f = \mathbb{Q}(\sqrt{d})$.

We choose $M = 22$ if the degree of the polynomial $P(x)$ is 5, and $M = 16$ if it is 6. We have designed two programs, one for the case of $f$ having an extra twist, and the other for $\varepsilon = 1$. Notice that the case such that $f$ has an extra twist and $\varepsilon = 1$ is included in both programs. These two programs return a set of possible solutions, some of which differ only in the term $\{a_n\}_{2 \leq n \leq M}$. To simplify the sieve process, we take the solutions of the form $(P(X), \{a_n\}_{2 \leq n \leq M''}, k, d)$ where $M'' = 11$ or $M'' = 7$, depending on whether the degree of the polynomial $P(x)$ is 5 or 6; if $(P(X), \{a_n\}_{2 \leq n \leq M''}, k, d)$ is a possible solution then we do not take in account the solution obtained by changing $\{a_n\}_{2 \leq n \leq M''}$ by $\{{}^\sigma a_n\}_{2 \leq n \leq M''}$ since $A_f = A_{\sigma f}$. Below, we give the number of possible solutions in two tables constructed according to the value of $d$, the degree of $P$ and the order of $\varepsilon$ when $f$ has an extra twist:

$\varepsilon = 1$

| $d$ | deg $P = 5$ | deg $P = 6$ |
|---|---|---|
| 2 | 369 | 17 |
| 3 | 302 | 7 |
| 5 | 386 | 28 |
| 6 | 22 | 0 |
| 7 | 15 | 0 |
| 10 | 0 | – |
| 11 | 0 | – |
| 13 | 116 | 3 |
| 17 | 98 | 1 |
| 21 | 44 | 0 |
| 29 | 8 | 0 |
| 33 | 0 | – |
| 37 | 0 | – |
| 41 | 0 | – |

Total: 1416

With an extra twist

| ord $\varepsilon$ | $d$ | deg $P = 5$ | deg $P = 6$ |
|---|---|---|---|
| 1 | 2 | 131 | 4 |
| 1 | 3 | 138 | 6 |
| 1 | 5 | 8 | 0 |
| 1 | 6 | 6 | 0 |
| 1 | 7 | 5 | 0 |
| 1 | 10 | 0 | – |
| 1 | 11 | 0 | – |
| 2 | $-1$ | 351 | 6 |
| 2 | $-2$ | 39 | 4 |
| 2 | $-3$ | 154 | 4 |
| 2 | $-5$ | 18 | 1 |
| 2 | $-6$ | 8 | 0 |
| 2 | $-7$ | 19 | 1 |
| 2 | $-10$ | 0 | – |
| 2 | $-11$ | 0 | – |
| 3 | $-3$ | 4 | 1 |
| 4 | $-1$ | 373 | 7 |
| 6 | $-3$ | 458 | 6 |

Total: 1752

## 4. Discarding

First we will sieve the possible solutions with an extra twist and next we will focus our attention to those without one.

4.1. **Extra twist.** We give six sieves that we will use to discard solutions returned by our program (with an extra twist) which do not satisfy Theorem 2.3.

*First Sieve*

The following proposition provides us with a criterion to make the first sieve.

**Proposition 4.1.** *Let $f \in S_2(N, \varepsilon)$ be a normalized newform with an extra twist and such that* dim $A_f = 2$. *Put*

$$n := \begin{cases} \text{ord}\, \varepsilon & \text{if} \quad \varepsilon \neq 1 \\ 2 & \text{if} \quad \varepsilon = 1 \end{cases}$$

*Then, for all primes $p$ that do not divide $N$, the characteristic polynomial $Q_{p^n}(t)$ of the Frobenius endomorphism $\mathrm{Frob}_{p^n}$ acting on the Tate module of $A_{f/\mathbb{F}_p}$ is of the form:*

(3) $$Q_{p^n}(t) = (t^2 - a\,t + p^n)^2, \quad \text{with} \quad a \in \mathbb{Z}.$$

*Proof.* This proposition is a particular case of Proposition 6.2 of [BG97] for the case dim $A_f = 2$ although here $f$ can have CM. By Eichler-Shimura, if $p \nmid N$ then there exists an algebraic integer $\alpha_p$ such that:

$$a_p = \alpha_p + \overline{\alpha_p}\varepsilon(p), \quad \alpha_p\overline{\alpha_p} = p, \quad Q_p(t) = (t^2 - a_p t + p\,\varepsilon(p))(t^2 - {}^\sigma a_p t + p^\sigma \varepsilon(p)).$$

Let $\chi$ be a character associated to an extra twist, that we take $\chi = \varepsilon^{-1}$ when $\varepsilon \neq 1$. It can be easily checked that:
$$a_p^n \in \mathbb{Z}, \quad a_p^2 \chi(p) \in \mathbb{Z}.$$

Put $a := \alpha_p^n + (\overline{\alpha_p}\varepsilon(p))^n$. Then $a = \alpha_p^n + (\overline{\alpha_p}\chi(p)^{-1})^n$ and we obtain (cf. Lemma 6.1 [BG97]) that:

$$a = \begin{cases} a_p^n \left( \sum_{i=1}^{[n/2]} (-1)^i \frac{n}{n-i} p^i (a_p^2 \chi(p))^{-i} \right) & \text{if} \quad a_p \neq 0, \\ 0 & \text{if} \quad a_p = 0 \quad \text{and} \quad n \equiv 1 \pmod 2, \\ 2(-1)^{n/2}(p\chi(p)^{-1})^{n/2} & \text{if} \quad a_p = 0 \quad \text{and} \quad n \equiv 0 \pmod 2. \end{cases}$$

Consequently, we obtain that $a \in \mathbb{Z}$ and $Q_{p^n}(t) = (t^2 - a\,t + p^n)^2$. $\square$

Applying the above proposition, we drop the curves obtained previously such that the corresponding polynomial $Q_{p^n}(t)$ for $J(C)$ is not of the form (3) for some prime $p < 100$ that does not divide the discriminant of the polynomial $P(x)$. After this sieve, the number of possible solutions has decreased from 1752 as far as 288.

### Second Sieve

This sieve is based on the fact (cf. [ST61]) that for every prime $p$ not dividing the conductor of $A_f$ over $\mathbb{Q}$, $\mathcal{N}_{\mathbb{Q}}(A_f)$, we have the following inclusion
$$\mathbb{Q} \otimes \mathrm{End}_{\mathbb{Q}}(A_f) \hookrightarrow \mathbb{Q} \otimes \mathrm{End}_{\mathbb{F}_p}(A_{f/\mathbb{F}_p}).$$

We have discarded the hyperelliptic curves such that the corresponding field $\mathbb{Q}(\sqrt{d})$ is not contained in $\mathbb{Q} \otimes \mathrm{End}_{\mathbb{F}_p}(J(C)/\mathbb{F}_p)$ for some $p \leq 29$ and $p$ not dividing the discriminant of $C$.

Now, only 94 possible solutions remain. After computing Clebsch's invariants (cf. [Bol88]), we have checked that every one of these correponds to a curve with a nonhyperelliptic involution and, thus, its jacobian is not $\overline{\mathbb{Q}}$-simple. Moreover, we have computed all their nonhyperelliptic involutions.

### Third Sieve

The following lemma will provide a criterion to eliminate some curves whose jacobians are not $\mathbb{Q}$-simple.

**Lemma 4.2** ([CGLR99]). *If $C$ is a curve defined over $\mathbb{Q}$ of genus 2 which has a nonhyperelliptic involution defined over a number field $K$, then*
$$J(C) \stackrel{K}{\sim} E \times E',$$
*where $E$ and $E'$ are elliptic curves defined over $K$.*

We discard the curves that admit a nonhyperelliptic involution defined over $\mathbb{Q}$, since their jacobians are not $\mathbb{Q}$-simple. After this sieve, 79 possible solutions remain.

### Fourth Sieve

Let $f \in S_2(N, \varepsilon)$ and $A_f$ be as above. Let $p$ be a prime that does not divide N and, as in Proposition 4.1, let us denote by $Q_p(t)$ the characteristic polynomial of the Frobenius endomorphism $\mathrm{Frob}_p$ acting on the Tate module of $A_{f/\mathbb{F}_p}$. By Eichler-Shimura, we know that

$$Q_p(t) = (t^2 - a_p\, t + \varepsilon(p)\, p)(t^2 - {}^\sigma a_p\, t + {}^\sigma\varepsilon(p)\, p)\,,$$

where $a_p$ denotes the corresponding coefficient in the $q$-expansion of $f$. If $p$ does not divide the discriminant of the polynomial $P(x)$, $a_p \neq 0$ and $\varepsilon \neq 1$, $\varepsilon(p)$ is determined by $a_p$, since $\varepsilon(p) = a_p/{}^\sigma a_p$. Recall that for each curve $C$, we have obtained the polynomial $P(x)$ and the values $\varepsilon(2), \varepsilon(3), a_2, \ldots, a_{M''}$. We discard the curves such that for some prime $p \leq M''$, $a_p \neq 0$, $p$ does not divide the discriminant of $P(x)$, and the polynomial $Q_p(t)$ attached to $J(C)$ does not coincide with the polynomial $(t^2 - a_p\, t + \varepsilon(p)\, p)(t^2 - {}^\sigma a_p\, t + {}^\sigma\varepsilon(p)\, p)$. By this stage, only 64 possible solutions have passed the previous four sieves.

### Fifth Sieve

The following theorem provides information about the conductor of the jacobian of $C$, when it is modular.

**Theorem 4.3** ([Car86]). *Let $f \in S_2(\Gamma_1(N))$ be a normalized newform such that $\dim A_f = d$, then*

$$\mathcal{N}_\mathbb{Q}(A_f) = N^d\,,$$

*where, as before, $\mathcal{N}_\mathbb{Q}(A_f)$ denotes the conductor of $A_f$ over $\mathbb{Q}$.*

In our case $A_f \stackrel{\mathbb{Q}}{\sim} J(C)$ and $\dim A_f = 2$, and we delete the curves $C$ such that $\mathcal{N}_\mathbb{Q}(J(C))$ is not a square integer. We have used this sieve to discard the curves such that the odd factor of their conductor is not a rational square. We have computed this odd factor using Quing Liu's program `genus2reduction`. This sieve has allowed us to delete 2 of those 64 possible solutions.

### Sixth Sieve

In order to prove that the curves that remain after passing all the previous sieves are solutions of Theorem 2.3, we will use a criterion that we describe in section 6. Notice that we have obtained pairs of possible solutions $(P(x), f, k, d)$, $(P'(x), f', k', d')$ such that the corresponding attached curves are $\mathbb{Q}$-isomorphic. The following result we will allow us to solve this.

**Proposition 4.4** ([Rib77]). *Let $f_1, f_2 \in S_2(\Gamma_1(N))$ be normalized newforms with Nebentypus $\varepsilon_1, \varepsilon_2$ respectively and such that $A_{f_1} \stackrel{\mathbb{Q}}{\sim} A_{f_2}$. Then there exists $\sigma \in \mathrm{Gal}(\overline{\mathbb{Q}}/\mathbb{Q})$ such that $f_2 = {}^\sigma f_1$ and $\varepsilon_2 = {}^\sigma \varepsilon_1$.*

In particular, when $(P(x), f, k, d)$ is a solution to our problem, we can discard the other solutions with $\mathbb{Q}$-isomorphic curves to $y^2 = P(x)$. After the six sieves we have obtained 53 equations.

4.2. **Without extra twist.** Now, we are only interested in the cases that $f$ does not have an extra twist and this condition implies that $A_f$ is $\overline{\mathbb{Q}}$-simple (cf. [Rib80], [Mom81]). So, we will restrict our attention to the solutions such that $J(C)$ is $\overline{\mathbb{Q}}$-simple.

The first step is to discard all the curves which have a nonhyperelliptic involution from among the equations obtained for the case $\varepsilon = 1$. To do so we compute Clebsch's invariants. Next, we use the second, the fourth and the sixth sieve, as in the previous subsection. We finally obtain 96 equations.

## 5. Searching

In this paragraph we explain the method used to find the newform attached to a primitive modular curve that has passed all the sieves. In order to do so, we determine the possible levels $N$ and Nebentypus $\varepsilon$ of $f$. Next, among the newforms in $S_2(N,\varepsilon)$ we look for those whose first Fourier coefficients agree with the values obtained by our program.

We know that $N = \sqrt{\mathcal{N}_{\mathbb{Q}}(J(C))}$. Let $N = 2^k M$ with $M \in \mathbb{Z}$ and $(2, M) = 1$. We compute the odd factor $M$ using Quing Liu's program `genus2reduction`. The algorithm described in [Liu94] also allows computation of the even factor $2^k$, but since this algorithm is relatively complicated we have used an alternative route. Proposition 5.5 of [Bru95] provides bounds for the values $k$. Indeed, let $m = \lceil \frac{k}{2} - 2 \rceil$ and $\zeta$ be a primitive $2^m$th root of the unity, then $\mathbb{Q}(\varepsilon, \zeta^2) \subseteq K_f$. Thus, for every quadratic field $K_f$ there are a few values of $k$ to search for our $f$ in the level $N = 2^k M$.

Among the curves we have found, we distinguish the following cases:

(i) *Case $\varepsilon = 1$ and without extra twist.* If $\deg P = 6$ then $2 \nmid N$ and, therefore, the conductor of $J(C)$ can be computed by the `genus2reduction` program. Otherwise, we use the Brumer's bound to compute the $2^k$-factor.

(ii) *Case ord $\varepsilon \le 2$ and with extra twist.* In this case, $C$ always has a nonhyperelliptic involution defined over a quadratic field $K$ and the computation of the conductor of $J(C)$ is as follows. We find an elliptic curve $E/K$ such that $J(C) \stackrel{K}{\sim} E^2$. Then, by the universal property of Weil's restriction of scalars, $J(C)$ is $\mathbb{Q}$-isogenous to $\mathrm{Res}_{K/\mathbb{Q}}(E)$ and its conductor is ([Mil72])

$$\mathcal{N}_{\mathbb{Q}}(J(C)) = N_{K/\mathbb{Q}}(\mathcal{N}_K(E)) \cdot d_{K/\mathbb{Q}}^2.$$

We note that if ord $\varepsilon = 2$ then the field $K$ is real and the Nebentypus $\varepsilon$ should be the character attached to the quadratic field $K$.

(iii) *Case ord $\varepsilon > 2$.* In this case, all the curves which have passed all the sieves are without CM and the field $L_C$ (the field of definition of all the involutions on $C$) must be equal to the field fixed by $\varepsilon$, $L$, since $[L_C : \mathbb{Q}] = \mathrm{ord}\,\varepsilon$. This fact determines $\varepsilon$ up to Galois conjugation. If $2^k$ is unknown then we use [Bru95] and the fact that $\mathfrak{f}_\varepsilon | 2^k M$, where $\mathfrak{f}_\varepsilon$ is the conductor of $\varepsilon$.

## 6. Checking

We consider the congruence subgroup:

$$\Gamma_\varepsilon(N) = \left\{ \begin{pmatrix} A & B \\ C & D \end{pmatrix} \in \Gamma_0(N) : \varepsilon(D) = 1 \right\},$$

and denote by $X(N, \varepsilon)$ the corresponding modular curve. We have that

$$S_2(\Gamma_\varepsilon(N)) = S_2(\Gamma_0(N)) \oplus \left( \bigoplus_{k=1}^{n-1} S_2(N, \varepsilon^k) \right),$$

where $n$ is the order of $\varepsilon$. Then, the genus of $X(N, \varepsilon)$ is

$$g = \begin{cases} g_0 + \sum_{k=1}^{n-1} \dim S_2(N, \varepsilon^k) & \text{if } \varepsilon \ne 1, \\ g_0 & \text{if } \varepsilon = 1, \end{cases}$$

where $g_0$ denotes the genus of $X_0(N)$.

Now, we assume that the modular form $f \in S_2(N, \varepsilon)$ is known. Let $h_1, h_2, x$ and $y$ be as in (1), (2) of the Section 2. Put

$$\mathrm{div}(h_2\, dq/q) = \sum_{i=1}^{h} m_i\, P_i, \quad \text{with} \quad P_i \in X_{\Gamma_\varepsilon} \quad \text{and} \quad P_i \neq P_j \quad \text{if} \quad i \neq j.$$

By the Riemann-Roch Theorem, we have $\deg(\mathrm{div}(h_2\, dq/q)) = \sum_{i=1}^{h} m_i = 2\,g - 2$. Let us denote by $\mathrm{div}^-(*)$ the positive divisor which is the polar part of $\mathrm{div}(*)$. Then

$$\begin{aligned}
\mathrm{div}^-(x) &= \mathrm{div}^-(h_1/h_2) \leq \mathrm{div}(h_2\, dq/q) \\
\mathrm{div}^-(y) &= \mathrm{div}^-(dx/h_2) \leq \mathrm{div}^-(dx) + \mathrm{div}(h_2\, dq/q) \\
&\leq \sum_{i=1}^{h}(m_i+1)P_i + \sum_{i=1}^{h} m_i P_i = 2\,\mathrm{div}(h_2) + \sum_{i=1}^{h} P_i.
\end{aligned}$$

Then, $\deg(\mathrm{div}^-(x)) \leq 2g-2$, $\deg(\mathrm{div}^-(y)) \leq 3(2g-2)$. Therefore,

$$\deg(\mathrm{div}^-(y^2 - P(x))) \leq 6(2g-2),$$

for every polynomial $P$ of degree $\leq 6$. Thus, we obtain the following criterion.

**Proposition 6.1.** *Keep the notations as above. If there exists a polynomial $P$ of degree $\leq 6$ such that*

$$y^2 - P(x) = O(q^{c_{N,\varepsilon}}) \quad \text{with} \quad c_{N,\varepsilon} = 6(2g-2) + 1,$$

*then $y^2 = P(x)$. If in addition $P$ is without double roots and has degree $5$ or $6$, then the hyperelliptic curve $C: y^2 = P(x)$ of genus two is primitive modular of level $N$.*

In order to obtain the first $c_{N,\varepsilon}$ coefficients of the Fourier expansion of the modular forms, we have used the program HECKE [Ste99] of W. A. Stein, a program designed by J. Quer, and some computations by M. Müller and W. A. Stein.

Finally, we obtain the following result which completes Theorem 2.3.

**Theorem 6.2.** *There are exactly $149$ primitive modular curves of genus two over $\mathbb{Q}$, up to $\mathbb{Q}$-isomorphism, such that their jacobians are $\mathbb{Q}$-simple.*

## 7. AN EXAMPLE

To conclude, we give an example of a genus 2 curve that is not primitive modular, but whose jacobian is.

Let us consider the hyperelliptic curve $C: y^2 = x(x-1)(x^3 - 3x + 1)$ and let $L$ be the number field $\mathbb{Q}(\lambda)$, where $\lambda$ is a root of the polynomial $x^3 - 3x + 1$. Since $C$ is a particular case of the curves treated in [CGLR99], we know that $J(C)$ is an abelian variety of $\mathrm{GL}_2$-type, $\mathbb{Q} \otimes \mathrm{End}_\mathbb{Q}(J(C)) = \mathbb{Q}(\sqrt{-3})$ and there exists a morphism $\pi: C \to E$ defined over $L$, where $E$ is the elliptic curve given by $y^2 = x^3 - 3^5\, 7\, x + 3^6\, 37$. We have that $\mathcal{N}_\mathbb{Q}(E) = 2^2\, 3^4$ and it is the curve 324A in Cremona's notation ([Cre97]). By the universal property of Weil's restriction of scalars we have

$$\mathrm{Res}_{L/\mathbb{Q}}(E) \stackrel{\mathbb{Q}}{\sim} J(C) \times E.$$

Now, let $\varepsilon$ be a Dirichlet character of order 3 and conductor 9. We take the normalized newform $f \in S_2(36, \varepsilon)$ with first Fourier coefficients:

$$f(q) = \sum_{n>0} a_n q^n = q - \sqrt{-3}\, q^3 + 3\, \frac{-1+\sqrt{-3}}{2}\, q^5 + \frac{1+\sqrt{-3}}{2}\, q^7 + \cdots.$$

The field fixed by the Nebentypus $\varepsilon$ of $f$ is the field $L$. It is easy to check that $\dim A_f = 2$ and $A_f$ is not $\overline{\mathbb{Q}}$-simple, i.e. $A_f$ is $L$-isogenous to the square of an elliptic curve defined over $L$. By Atkin-Li [AL78] we have that $g = f \otimes \varepsilon = \sum_{n>0} \varepsilon(n) a_n q^n$ is a newform of level $2^2\, 3^4$ with trivial Nebentypus and has Fourier $q$-expansion with integer coefficients. Then, the abelian variety $A_g$ attached by Shimura to $g$ is an elliptic curve defined over $\mathbb{Q}$ and we check that $A_g$ is $\mathbb{Q}$-isogenous to $E$. By ([Shi71]) we have that there exists a nonconstant morphism $\pi': A_f \to A_g$ defined over $L$. Thus, with the same reasoning as above we obtain

$$\mathrm{Res}_{L/\mathbb{Q}}(E) \stackrel{\mathbb{Q}}{\sim} A_f \times E\,,$$

and, therefore, $J(C) \stackrel{\mathbb{Q}}{\sim} A_f$. Thus, we have proved that $J(C)$ is $\mathbb{Q}$-simple primitive modular of level $N = 36$, but $C$ is not a primitive modular curve since it does not appear in the tables of the appendix.

## Appendix

Here we present three tables that are available to be taken electronically in the address `http://mat.uab.es/enrikegj`. We give the 149 primitive modular curves we have obtained in the first or in the second table, depending on whether the corresponding newform has an extra twist or not. These two tables have two columns. The first column contains the hyperelliptic equations. The levels $N$ and the Nebentypus $\varepsilon$ of the corresponding newforms are in the second column, where the superscript CM denotes complex multiplication. When there is a unique curve for a given level $N$, this is denoted by $C_N$; otherwise we use the notation $C_{N,A}, C_{N,B}, \ldots$. To give the character $\varepsilon$, if it is nontrivial, we enter a comma separated list of the degrees of the images of the factors corresponding to the primes which divide the level. There is some subtlety at 2 because $(\mathbb{Z}/2^n\mathbb{Z})^*$ is not necessarily cyclic. If $2|N$ but 4 does not divide $N$, ignore 2. Otherwise, give either 1 or 2 degrees for 2, if 4 or 8 divides $N$, respectively. This notation identify characters up to Galois conjugation. If the Nebentypus is trivial we use $S_2(N)$ to denote $S_2(N, 1)$.

Notice that the 149 curves obtained provide 99 isomorphism classes (over $\overline{\mathbb{Q}}$). We remark that when $d \equiv 1 \pmod 4$ we have made a change in the variables $x$, $y$ to obtain models defined over $\mathbb{Z}$. In this case, we have taken $h_1 - h_2, -2h_2$ instead $h_1, h_2$ to compute $x$ and $y$.

In the third table we give the coefficients $a_p$, $p \leq 13$, of the $q$-expansion of the newforms. For simplicity, $w_d$ denotes a root of $x^2 - x - (d-1)/4$ when $d \equiv 1 \pmod 4$.

Table 1: Non $\overline{\mathbb{Q}}$-simple Jacobian

| $C$ | : $y^2 = P(x)$ | $f$ | $\in$ | $S_2(N, \varepsilon)$ |
|---|---|---|---|---|
| $C_{13}$ | : $y^2 = x^6 + 4\,x^5 + 6\,x^4 + 2\,x^3 + x^2 + 2\,x + 1$ | $f_{13}$ | $\in$ | $S_2(13, [6])$ |
| $C_{16}$ | : $y^2 = x^6 + 2\,x^5 - x^4 - x^2 - 2\,x + 1$ | $f_{16}$ | $\in$ | $S_2(16, [1, 4])$ |
| $C_{18}$ | : $y^2 = x^6 + 2\,x^5 + 5\,x^4 + 10\,x^3 + 10\,x^2 + 4\,x + 1$ | $f_{18}$ | $\in$ | $S_2(18, [1, 3])$ |
| $C_{28}$ | : $y^2 = x^5 - 4\,x^4 - 13\,x^3 - 9\,x^2 - x$ | $f_{28}$ | $\in$ | $S_2(28, [1, 3])$ |
| $C_{40}$ | : $y^2 = x^5 - 3\,x^4 - 12\,x^2 - 16\,x$ | $f_{40}$ | $\in$ | $S_2(40, [1, 1, 2])$ |

| $C$ | : | $y^2 = P(x)$ | $f$ | $\in$ | $S_2(N, \varepsilon)$ |
|---|---|---|---|---|---|
| $C_{45}$ | : | $y^2 = x^6 + 22\,x^3 + 125$ | $f_{45}^{\mathrm{CM}}$ | $\in$ | $S_2(45, [1, 2])$ |
| $C_{48}$ | : | $y^2 = x^5 - 5\,x^4 + 20\,x^2 - 16\,x$ | $f_{48}^{\mathrm{CM}}$ | $\in$ | $S_2(48, [2, 1, 2])$ |
| $C_{52}$ | : | $y^2 = x^5 - 5\,x^3 - 5\,x^2 - x$ | $f_{52}$ | $\in$ | $S_2(52, [1, 6])$ |
| $C_{63}$ | : | $y^2 = x^6 - 26\,x^3 - 27$ | $f_{63}$ | $\in$ | $S_2(63)$ |
| $C_{64,A}$ | : | $y^2 = x^5 - 16\,x$ | $f_{64,A}^{\mathrm{CM}}$ | $\in$ | $S_2(64, [1, 2])$ |
| $C_{64,B}$ | : | $y^2 = x^5 - 2\,x^4 - 2\,x^2 - x$ | $f_{64,B}$ | $\in$ | $S_2(64, [1, 4])$ |
| $C_{80}$ | : | $y^2 = x^5 + 3\,x^4 + 12\,x^2 - 16\,x$ | $f_{80}$ | $\in$ | $S_2(80, [1, 1, 2])$ |
| $C_{81}$ | : | $y^2 = x^6 - 18\,x^3 - 27$ | $f_{81}$ | $\in$ | $S_2(81)$ |
| $C_{100}$ | : | $y^2 = x^5 + 5\,x^3 + 5\,x - 11$ | $f_{100}$ | $\in$ | $S_2(100, [1, 2])$ |
| $C_{112,A}$ | : | $y^2 = x^5 + 3\,x^4 + x^3 - 2\,x^2 - x$ | $f_{112,A}$ | $\in$ | $S_2(112, [2, 1, 6])$ |
| $C_{112,B}$ | : | $y^2 = x^5 - 3\,x^4 + x^3 + 2\,x^2 - x$ | $f_{112,B}$ | $\in$ | $S_2(112, [2, 1, 6])$ |
| $C_{112,C}$ | : | $y^2 = x^5 + 4\,x^4 - 13\,x^3 + 9\,x^2 - x$ | $f_{112,C}$ | $\in$ | $S_2(112, [1, 1, 3])$ |
| $C_{117}$ | : | $y^2 = x^6 - 10\,x^3 - 27$ | $f_{117}$ | $\in$ | $S_2(117)$ |
| $C_{128}$ | : | $y^2 = x^5 + 64\,x$ | $f_{128}^{\mathrm{CM}}$ | $\in$ | $S_2(128, [1, 2])$ |
| $C_{148}$ | : | $y^2 = x^5 + 8\,x^4 + 11\,x^3 + 3\,x^2 - x$ | $f_{148}$ | $\in$ | $S_2(148, [1, 6])$ |
| $C_{160,A}$ | : | $y^2 = x^5 - 12\,x^3 - 64\,x$ | $f_{160,A}$ | $\in$ | $S_2(160)$ |
| $C_{160,B}$ | : | $y^2 = x^5 - 11\,x^4 - 11\,x^2 - x$ | $f_{160,B}$ | $\in$ | $S_2(160, [2, 1, 4])$ |
| $C_{160,C}$ | : | $y^2 = x^5 + x^4 + x^2 - x$ | $f_{160,C}$ | $\in$ | $S_2(160, [2, 1, 4])$ |
| $C_{160,D}$ | : | $y^2 = x^5 - x^4 - x^2 - x$ | $f_{160,D}$ | $\in$ | $S_2(160, [2, 1, 4])$ |
| $C_{160,E}$ | : | $y^2 = x^5 + 11\,x^4 + 11\,x^2 - x$ | $f_{160,E}$ | $\in$ | $S_2(160, [2, 1, 4])$ |
| $C_{189}$ | : | $y^2 = x^6 - 2\,x^3 - 27$ | $f_{189}$ | $\in$ | $S_2(189)$ |
| $C_{192}$ | : | $y^2 = x^5 - 5\,x^4 + 20\,x^3 - 40\,x^2 + 44\,x - 20$ | $f_{192}^{\mathrm{CM}}$ | $\in$ | $S_2(192, [2, 1, 2])$ |
| $C_{208,A}$ | : | $y^2 = x^5 - 5\,x^3 + 5\,x^2 - x$ | $f_{208,A}$ | $\in$ | $S_2(208, [1, 1, 6])$ |
| $C_{243}$ | : | $y^2 = x^6 + 6\,x^3 - 27$ | $f_{243}$ | $\in$ | $S_2(243)$ |
| $C_{256}$ | : | $y^2 = x^5 - 64\,x$ | $f_{256}^{\mathrm{CM}}$ | $\in$ | $S_2(256)$ |
| $C_{320}$ | : | $y^2 = x^5 + 12\,x^3 - 64\,x$ | $f_{320}$ | $\in$ | $S_2(320)$ |
| $C_{400}$ | : | $y^2 = x^5 + 5\,x^3 + 5\,x + 11$ | $f_{400}$ | $\in$ | $S_2(400, [1, 1, 2])$ |
| $C_{512,A}$ | : | $y^2 = x^5 + 4\,x^3 - 4\,x$ | $f_{512,A}$ | $\in$ | $S_2(512)$ |
| $C_{512,B}$ | : | $y^2 = x^5 - 4\,x^3 - 4\,x$ | $f_{512,B}$ | $\in$ | $S_2(512)$ |
| $C_{544}$ | : | $y^2 = x^5 - x^3 - 4\,x$ | $f_{544}$ | $\in$ | $S_2(544)$ |
| $C_{592}$ | : | $y^2 = x^5 - 8\,x^4 + 11\,x^3 - 3\,x^2 - x$ | $f_{592}$ | $\in$ | $S_2(592, [1, 1, 6])$ |
| $C_{768,A}$ | : | $y^2 = x^5 - x^4 + 2\,x^3 + 2\,x^2 - 2\,x - 2$ | $f_{768,A}$ | $\in$ | $S_2(768, [2, 1, 2])$ |
| $C_{768,B}$ | : | $y^2 = x^5 + x^4 + 2\,x^3 - 2\,x^2 - 2\,x + 2$ | $f_{768,B}$ | $\in$ | $S_2(768, [2, 1, 2])$ |
| $C_{928,A}$ | : | $y^2 = x^5 + 5\,x^4 + 5\,x^2 - x$ | $f_{928,A}$ | $\in$ | $S_2(928, [2, 1, 4])$ |
| $C_{928,B}$ | : | $y^2 = x^5 - 5\,x^4 - 5\,x^2 - x$ | $f_{928,B}$ | $\in$ | $S_2(928, [2, 1, 4])$ |
| $C_{1088}$ | : | $y^2 = x^5 + x^3 - 4\,x$ | $f_{1088}$ | $\in$ | $S_2(1088)$ |
| $C_{1280,A}$ | : | $y^2 = x^5 + 2\,x^3 - 4\,x$ | $f_{1280,A}$ | $\in$ | $S_2(1280)$ |
| $C_{1280,B}$ | : | $y^2 = x^5 - 2\,x^3 - 4\,x$ | $f_{1280,B}$ | $\in$ | $S_2(1280)$ |

| $C$ | : | $y^2 = P(x)$ | $f$ | $\in$ | $S_2(N, \varepsilon)$ |
|---|---|---|---|---|---|
| $C_{1280,C}$ | : | $y^2 = x^5 + 22\,x^3 - 4\,x$ | $f_{1280,C}$ | $\in$ | $S_2(1280)$ |
| $C_{1280,D}$ | : | $y^2 = x^5 - 22\,x^3 - 4\,x$ | $f_{1280,D}$ | $\in$ | $S_2(1280)$ |
| $C_{1280,E}$ | : | $y^2 = x^5 + 8\,x^3 - 4\,x$ | $f_{1280,E}$ | $\in$ | $S_2(1280)$ |
| $C_{1280,F}$ | : | $y^2 = x^5 - 8\,x^3 - 4\,x$ | $f_{1280,F}$ | $\in$ | $S_2(1280)$ |
| $C_{1312}$ | : | $y^2 = x^5 - 5\,x^3 - 4\,x$ | $f_{1312}$ | $\in$ | $S_2(1312)$ |
| $C_{2080}$ | : | $y^2 = x^5 + 7\,x^3 - 4\,x$ | $f_{2080}$ | $\in$ | $S_2(2080)$ |
| $C_{2624}$ | : | $y^2 = x^5 + 5\,x^3 - 4\,x$ | $f_{2624}$ | $\in$ | $S_2(2624)$ |
| $C_{4160}$ | : | $y^2 = x^5 - 7\,x^3 - 4\,x$ | $f_{4160}$ | $\in$ | $S_2(4160)$ |
| $C_{7424,A}$ | : | $y^2 = x^5 + 10\,x^3 - 4\,x$ | $f_{7424,A}$ | $\in$ | $S_2(7424)$ |
| $C_{7424,B}$ | : | $y^2 = x^5 - 10\,x^3 - 4\,x$ | $f_{7424,B}$ | $\in$ | $S_2(7424)$ |

Table 2: $\overline{\mathbb{Q}}$-simple Jacobian

| $C$ | : | $y^2 = P(x)$ | $f$ | $\in$ | $S_2(N)$ |
|---|---|---|---|---|---|
| $C_{23}$ | : | $y^2 = x^6 - 8\,x^5 + 2\,x^4 + 2\,x^3 - 11\,x^2 + 10\,x - 7$ | $f_{23}$ | $\in$ | $S_2(23)$ |
| $C_{29}$ | : | $y^2 = x^6 + 2\,x^5 - 17\,x^4 - 66\,x^3 - 83\,x^2 - 32\,x - 4$ | $f_{29}$ | $\in$ | $S_2(29)$ |
| $C_{31}$ | : | $y^2 = x^6 - 14\,x^5 + 61\,x^4 - 106\,x^3 + 66\,x^2 - 8\,x - 3$ | $f_{31}$ | $\in$ | $S_2(31)$ |
| $C_{35}$ | : | $y^2 = x^6 - 4\,x^5 + 2\,x^4 - 32\,x^3 - 27\,x^2 - 64\,x - 76$ | $f_{35}$ | $\in$ | $S_2(35)$ |
| $C_{39}$ | : | $y^2 = x^6 + 6\,x^5 - 5\,x^4 - 66\,x^3 - 59\,x^2 - 12\,x$ | $f_{39}$ | $\in$ | $S_2(39)$ |
| $C_{67}$ | : | $y^2 = x^6 + 2\,x^5 + x^4 - 2\,x^3 + 2\,x^2 - 4\,x + 1$ | $f_{67}$ | $\in$ | $S_2(67)$ |
| $C_{68}$ | : | $y^2 = x^5 - 11\,x^4 + 7\,x^3 + 7\,x^2 - 12\,x + 8$ | $f_{68}$ | $\in$ | $S_2(68)$ |
| $C_{73}$ | : | $y^2 = x^6 + 2\,x^5 + x^4 + 6\,x^3 + 2\,x^2 - 4\,x + 1$ | $f_{73}$ | $\in$ | $S_2(73)$ |
| $C_{85}$ | : | $y^2 = x^6 + 2\,x^5 + 7\,x^4 + 6\,x^3 + 13\,x^2 - 8\,x + 4$ | $f_{85}$ | $\in$ | $S_2(85)$ |
| $C_{87}$ | : | $y^2 = x^6 - 6\,x^5 + 13\,x^4 - 18\,x^3 + 10\,x^2 - 3$ | $f_{87}$ | $\in$ | $S_2(87)$ |
| $C_{88}$ | : | $y^2 = x^5 - 8\,x^4 - 4\,x^3 + 36\,x^2 - 32\,x$ | $f_{88}$ | $\in$ | $S_2(88)$ |
| $C_{93}$ | : | $y^2 = x^6 + 6\,x^5 + 5\,x^4 - 6\,x^3 + 2\,x^2 + 1$ | $f_{93}$ | $\in$ | $S_2(93)$ |
| $C_{103}$ | : | $y^2 = x^6 + 6\,x^5 + 5\,x^4 + 2\,x^3 + 2\,x^2 + 1$ | $f_{103}$ | $\in$ | $S_2(103)$ |
| $C_{104}$ | : | $y^2 = x^5 - 8\,x^4 + x^3 + 30\,x^2 - 20\,x + 8$ | $f_{104}$ | $\in$ | $S_2(104)$ |
| $C_{107}$ | : | $y^2 = x^6 - 4\,x^5 + 10\,x^4 - 18\,x^3 + 17\,x^2 - 10\,x + 1$ | $f_{107}$ | $\in$ | $S_2(107)$ |
| $C_{115}$ | : | $y^2 = x^6 + 6\,x^5 + 5\,x^4 + 10\,x^3 + 2\,x^2 + 1$ | $f_{115}$ | $\in$ | $S_2(115)$ |
| $C_{125}$ | : | $y^2 = x^6 - 4\,x^5 + 10\,x^4 - 10\,x^3 + 5\,x^2 + 2\,x - 3$ | $f_{125}$ | $\in$ | $S_2(125)$ |
| $C_{133}$ | : | $y^2 = x^6 + 10\,x^5 + 17\,x^4 + 14\,x^3 + 10\,x^2 + 4\,x + 1$ | $f_{133}$ | $\in$ | $S_2(133)$ |
| $C_{135}$ | : | $y^2 = x^6 - 6\,x^5 + 21\,x^4 - 54\,x^3 + 90\,x^2 - 108\,x + 45$ | $f_{135}$ | $\in$ | $S_2(135)$ |
| $C_{136}$ | : | $y^2 = x^5 - 19\,x^3 - 14\,x^2 + 28\,x - 8$ | $f_{136}$ | $\in$ | $S_2(136)$ |
| $C_{147}$ | : | $y^2 = x^6 + 6\,x^5 + 11\,x^4 + 6\,x^3 + 5\,x^2 + 4$ | $f_{147}$ | $\in$ | $S_2(147)$ |
| $C_{161}$ | : | $y^2 = x^6 + 2\,x^4 - 6\,x^3 + 17\,x^2 - 18\,x + 5$ | $f_{161}$ | $\in$ | $S_2(161)$ |
| $C_{165}$ | : | $y^2 = x^6 + 6\,x^5 + 11\,x^4 + 14\,x^3 + 5\,x^2 - 12\,x$ | $f_{165}$ | $\in$ | $S_2(165)$ |

| $C$ | : | $y^2 = P(x)$ | $f$ | $\in$ | $S_2(N)$ |
|---|---|---|---|---|---|
| $C_{167}$ | : | $y^2 = x^6 - 4x^5 + 2x^4 - 2x^3 - 3x^2 + 2x - 3$ | $f_{167}$ | $\in$ | $S_2(167)$ |
| $C_{175}$ | : | $y^2 = x^6 + 2x^5 - 3x^4 + 6x^3 - 14x^2 + 8x - 3$ | $f_{175}$ | $\in$ | $S_2(175)$ |
| $C_{176}$ | : | $y^2 = x^5 + 3x^4 - 26x^3 + 14x^2 + x + 7$ | $f_{176}$ | $\in$ | $S_2(176)$ |
| $C_{177}$ | : | $y^2 = x^6 + 2x^4 - 6x^3 + 5x^2 - 6x + 1$ | $f_{177}$ | $\in$ | $S_2(177)$ |
| $C_{184}$ | : | $y^2 = x^5 - 14x^3 - 7x^2 + 19x - 7$ | $f_{184}$ | $\in$ | $S_2(184)$ |
| $C_{188,A}$ | : | $y^2 = x^5 - x^4 + x^3 + x^2 - 2x + 1$ | $f_{188,A}$ | $\in$ | $S_2(188)$ |
| $C_{188,B}$ | : | $y^2 = x^5 - 5x^4 + 5x^3 - 15x^2 + 6x - 11$ | $f_{188,B}$ | $\in$ | $S_2(188)$ |
| $C_{191}$ | : | $y^2 = x^6 + 2x^4 + 2x^3 + 5x^2 - 6x + 1$ | $f_{191}$ | $\in$ | $S_2(191)$ |
| $C_{205}$ | : | $y^2 = x^6 + 2x^4 + 10x^3 + 5x^2 - 6x + 1$ | $f_{205}$ | $\in$ | $S_2(205)$ |
| $C_{207}$ | : | $y^2 = x^6 + 6x^5 + 3x^4 - 26x^3 - 27x^2 - 12x$ | $f_{207}$ | $\in$ | $S_2(207)$ |
| $C_{208,B}$ | : | $y^2 = x^5 + 3x^4 - 21x^3 + 5x^2 + 16x - 12$ | $f_{208,B}$ | $\in$ | $S_2(208)$ |
| $C_{209}$ | : | $y^2 = x^6 - 4x^5 + 8x^4 - 8x^3 + 8x^2 + 4x + 4$ | $f_{209}$ | $\in$ | $S_2(209)$ |
| $C_{213}$ | : | $y^2 = x^6 + 2x^4 + 2x^3 - 7x^2 + 6x - 3$ | $f_{213}$ | $\in$ | $S_2(213)$ |
| $C_{221}$ | : | $y^2 = x^6 + 4x^5 + 2x^4 + 6x^3 + x^2 - 2x + 1$ | $f_{221}$ | $\in$ | $S_2(221)$ |
| $C_{224,A}$ | : | $y^2 = x^5 - 13x^4 + 48x^3 - 36x^2 - 32x + 32$ | $f_{224,A}$ | $\in$ | $S_2(224)$ |
| $C_{224,B}$ | : | $y^2 = x^5 + 3x^4 - 16x^3 - 20x^2 + 64x - 32$ | $f_{224,B}$ | $\in$ | $S_2(224)$ |
| $C_{261}$ | : | $y^2 = x^6 - 6x^4 + 10x^3 + 21x^2 - 30x + 9$ | $f_{261}$ | $\in$ | $S_2(261)$ |
| $C_{272,A}$ | : | $y^2 = x^5 - 10x^4 + 21x^3 + 48x^2 - 176x + 128$ | $f_{272,A}$ | $\in$ | $S_2(272)$ |
| $C_{272,B}$ | : | $y^2 = x^5 + 11x^4 + 7x^3 - 7x^2 - 12x - 8$ | $f_{272,B}$ | $\in$ | $S_2(272)$ |
| $C_{275}$ | : | $y^2 = x^6 - 4x^5 + 2x^4 - 22x^3 - 15x^2 - 30x - 35$ | $f_{275}$ | $\in$ | $S_2(275)$ |
| $C_{280}$ | : | $y^2 = x^5 - 5x^4 + 3x^3 + 9x^2 - 20x$ | $f_{280}$ | $\in$ | $S_2(280)$ |
| $C_{287}$ | : | $y^2 = x^6 - 4x^5 + 2x^4 + 6x^3 - 15x^2 + 14x - 7$ | $f_{287}$ | $\in$ | $S_2(287)$ |
| $C_{297}$ | : | $y^2 = x^6 + 6x^5 + 3x^4 - 36x^3 - 69x^2 - 54x - 15$ | $f_{297}$ | $\in$ | $S_2(297)$ |
| $C_{299}$ | : | $y^2 = x^6 - 10x^5 + 41x^4 - 78x^3 + 66x^2 - 28x + 5$ | $f_{299}$ | $\in$ | $S_2(299)$ |
| $C_{315}$ | : | $y^2 = x^6 + 6x^5 + 3x^4 - 18x^3 - 27x^2 - 24x - 4$ | $f_{315}$ | $\in$ | $S_2(315)$ |
| $C_{351}$ | : | $y^2 = x^6 - 6x^4 + 18x^3 + 9x^2 - 18x + 5$ | $f_{351}$ | $\in$ | $S_2(351)$ |
| $C_{357}$ | : | $y^2 = x^6 + 8x^4 - 8x^3 + 20x^2 - 12x + 12$ | $f_{357}$ | $\in$ | $S_2(357)$ |
| $C_{368}$ | : | $y^2 = x^5 - 5x^4 - 4x^3 + 39x^2 - 32x + 8$ | $f_{368}$ | $\in$ | $S_2(368)$ |
| $C_{376,A}$ | : | $y^2 = x^5 - 6x^4 + 11x^3 - 4x^2 - 2x + 1$ | $f_{376,A}$ | $\in$ | $S_2(376)$ |
| $C_{376,B}$ | : | $y^2 = x^5 - x^3 + 2x^2 - 2x + 1$ | $f_{376,B}$ | $\in$ | $S_2(376)$ |
| $C_{380}$ | : | $y^2 = x^5 - 5x^4 + 3x^3 + 7x^2 - 3x + 2$ | $f_{380}$ | $\in$ | $S_2(380)$ |
| $C_{416,A}$ | : | $y^2 = x^5 - 3x^4 + 7x^3 - x^2 - 8x + 20$ | $f_{416,A}$ | $\in$ | $S_2(416)$ |
| $C_{416,B}$ | : | $y^2 = x^5 - 2x^4 + 5x^3 - 12x^2 + 4x - 16$ | $f_{416,B}$ | $\in$ | $S_2(416)$ |
| $C_{440,A}$ | : | $y^2 = x^5 - 2x^3 - 7x^2 - 8x + 8$ | $f_{440,A}$ | $\in$ | $S_2(440)$ |
| $C_{440,B}$ | : | $y^2 = x^5 - 5x^4 + 12x^3 - 27x^2 + 25x - 30$ | $f_{440,B}$ | $\in$ | $S_2(440)$ |
| $C_{448,A}$ | : | $y^2 = x^5 - 7x^4 + 28x^3 - 72x^2 + 76x - 12$ | $f_{448,A}$ | $\in$ | $S_2(448)$ |
| $C_{448,B}$ | : | $y^2 = x^5 - 3x^4 + 12x^3 - 8x^2 - 20x + 4$ | $f_{448,B}$ | $\in$ | $S_2(448)$ |
| $C_{476,A}$ | : | $y^2 = x^5 - 3x^4 + 5x^3 - x^2 - 2x + 1$ | $f_{476,A}$ | $\in$ | $S_2(476)$ |

| $C$ | : | $y^2 = P(x)$ | $f$ | $\in$ | $S_2(N)$ |
|---|---|---|---|---|---|
| $C_{476,B}$ | : | $y^2 = x^5 - 2x^4 + 3x^3 - 6x^2 - 7$ | $f_{476,B}$ | $\in$ | $S_2(476)$ |
| $C_{525}$ | : | $y^2 = x^6 + 2x^4 - 10x^3 - 7x^2 - 30x + 9$ | $f_{525}$ | $\in$ | $S_2(525)$ |
| $C_{560}$ | : | $y^2 = x^5 - 7x^3 + 2x^2 - 8x + 12$ | $f_{560}$ | $\in$ | $S_2(560)$ |
| $C_{621}$ | : | $y^2 = x^6 - 6x^5 + 21x^4 - 42x^3 + 42x^2 - 24x + 5$ | $f_{621}$ | $\in$ | $S_2(621)$ |
| $C_{640,A}$ | : | $y^2 = x^5 - 10x^4 + 42x^3 - 96x^2 + 112x - 64$ | $f_{640,A}$ | $\in$ | $S_2(640)$ |
| $C_{640,B}$ | : | $y^2 = x^5 + 2x^3 + 4x^2 - 8x + 16$ | $f_{640,B}$ | $\in$ | $S_2(640)$ |
| $C_{645}$ | : | $y^2 = x^6 + 8x^4 + 20x^2 + 12x + 4$ | $f_{645}$ | $\in$ | $S_2(645)$ |
| $C_{704,A}$ | : | $y^2 = x^5 - 2x^4 + 10x^3 - 18x^2 + 16x - 16$ | $f_{704,A}$ | $\in$ | $S_2(704)$ |
| $C_{704,B}$ | : | $y^2 = x^5 - 3x^4 + 12x^3 - 10x^2 + 7x + 9$ | $f_{704,B}$ | $\in$ | $S_2(704)$ |
| $C_{752,A}$ | : | $y^2 = x^5 - 5x^3 + 20x^2 - 24x + 19$ | $f_{752,A}$ | $\in$ | $S_2(752)$ |
| $C_{752,B}$ | : | $y^2 = x^5 - 4x^4 + 7x^3 - 8x^2 + 4x - 1$ | $f_{752,B}$ | $\in$ | $S_2(752)$ |
| $C_{752,C}$ | : | $y^2 = x^5 + x^4 - 3x^3 - 3x^2 + 4x - 1$ | $f_{752,C}$ | $\in$ | $S_2(752)$ |
| $C_{752,D}$ | : | $y^2 = x^5 - 5x^4 + 9x^3 - 9x^2 + 4x - 1$ | $f_{752,D}$ | $\in$ | $S_2(752)$ |
| $C_{783}$ | : | $y^2 = x^6 - 6x^4 + 10x^3 - 15x^2 + 6x - 3$ | $f_{783}$ | $\in$ | $S_2(783)$ |
| $C_{880,A}$ | : | $y^2 = x^5 - 5x^4 + 8x^3 + 3x^2 - 23x + 8$ | $f_{880,A}$ | $\in$ | $S_2(880)$ |
| $C_{880,B}$ | : | $y^2 = x^5 + 2x^3 + 11x^2 - 8x + 24$ | $f_{880,B}$ | $\in$ | $S_2(880)$ |
| $C_{952}$ | : | $y^2 = x^5 - 8x^4 + 27x^3 - 30x^2 + 18x - 7$ | $f_{952}$ | $\in$ | $S_2(952)$ |
| $C_{1053}$ | : | $y^2 = x^6 - 6x^4 + 18x^3 - 27x^2 + 18x - 7$ | $f_{1053}$ | $\in$ | $S_2(1053)$ |
| $C_{1120,A}$ | : | $y^2 = x^5 + 5x^4 + 17x^3 + 15x^2 - 14x$ | $f_{1120,A}$ | $\in$ | $S_2(1120)$ |
| $C_{1120,B}$ | : | $y^2 = x^5 - 5x^4 + 17x^3 - 15x^2 - 14x$ | $f_{1120,B}$ | $\in$ | $S_2(1120)$ |
| $C_{1280,G}$ | : | $y^2 = x^5 + 5x^4 + 22x^3 + 14x^2 - 3x + 1$ | $f_{1280,G}$ | $\in$ | $S_2(1280)$ |
| $C_{1280,H}$ | : | $y^2 = x^5 - 5x^4 + 22x^3 - 14x^2 - 3x - 1$ | $f_{1280,H}$ | $\in$ | $S_2(1280)$ |
| $C_{1520}$ | : | $y^2 = x^5 + 5x^4 + 3x^3 - 7x^2 - 3x - 2$ | $f_{1520}$ | $\in$ | $S_2(1520)$ |
| $C_{1792,A}$ | : | $y^2 = x^5 - 3x^4 - 2x^3 + 2x^2 + 10x + 6$ | $f_{1792,A}$ | $\in$ | $S_2(1792)$ |
| $C_{1792,B}$ | : | $y^2 = x^5 + 3x^4 - 2x^3 - 2x^2 + 10x - 6$ | $f_{1792,B}$ | $\in$ | $S_2(1792)$ |
| $C_{1820}$ | : | $y^2 = x^5 + 7x^3 - 16x^2 + 12x + 1$ | $f_{1820}$ | $\in$ | $S_2(1820)$ |
| $C_{1904,A}$ | : | $y^2 = x^5 - 2x^4 + 3x^3 - 6x^2 + 4x - 1$ | $f_{1904,A}$ | $\in$ | $S_2(1904)$ |
| $C_{1904,B}$ | : | $y^2 = x^5 + 2x^4 + 3x^3 + 6x^2 + 7$ | $f_{1904,B}$ | $\in$ | $S_2(1904)$ |
| $C_{1904,C}$ | : | $y^2 = x^5 + 3x^4 + 5x^3 - 13x^2 + 12x - 1$ | $f_{1904,C}$ | $\in$ | $S_2(1904)$ |
| $C_{1916}$ | : | $y^2 = x^5 + 3x^4 - 7x^3 + 5x^2 - 2x + 1$ | $f_{1916}$ | $\in$ | $S_2(1916)$ |
| $C_{2240,A}$ | : | $y^2 = x^5 - 5x^4 + 17x^3 - 33x^2 + 40x - 30$ | $f_{2240,A}$ | $\in$ | $S_2(2240)$ |
| $C_{2240,B}$ | : | $y^2 = x^5 + 5x^4 + 17x^3 + 33x^2 + 40x + 30$ | $f_{2240,B}$ | $\in$ | $S_2(2240)$ |
| $C_{3159}$ | : | $y^2 = x^6 - 6x^5 + 21x^4 - 26x^3 + 18x^2 - 3$ | $f_{3159}$ | $\in$ | $S_2(3159)$ |
| $C_{7280}$ | : | $y^2 = x^5 - 5x^4 + 17x^3 - 15x^2 + 6x - 5$ | $f_{7280}$ | $\in$ | $S_2(7280)$ |
| $C_{7664}$ | : | $y^2 = x^5 - 3x^4 - 7x^3 - 5x^2 - 2x - 1$ | $f_{7664}$ | $\in$ | $S_2(7664)$ |

Table 3: First coefficients of the newforms

| $f$ | $a_2$ | $a_3$ | $a_5$ | $a_7$ | $a_{11}$ | $a_{13}$ |
|---|---|---|---|---|---|---|
| $f_{13}$ | $-w_{-3}-1$ | $2w_{-3}-2$ | $-2w_{-3}+1$ | $0$ | $0$ | $-3w_{-3}-1$ |
| $f_{16}$ | $-i-1$ | $i-1$ | $-i-1$ | $-2i$ | $i+1$ | $i-1$ |
| $f_{18}$ | $-w_{-3}$ | $w_{-3}-2$ | $0$ | $-2w_{-3}$ | $3w_{-3}$ | $2w_{-3}-2$ |
| $f_{23}$ | $w_5-1$ | $-2w_5+1$ | $2w_5-2$ | $2w_5$ | $-2w_5-2$ | $3$ |
| $f_{28}$ | $0$ | $-w_{-3}$ | $3w_{-3}-3$ | $-2w_{-3}-1$ | $3w_{-3}$ | $2$ |
| $f_{29}$ | $\sqrt{2}-1$ | $-\sqrt{2}+1$ | $-1$ | $2\sqrt{2}$ | $\sqrt{2}+1$ | $2\sqrt{2}-1$ |
| $f_{31}$ | $w_5$ | $-2w_5$ | $1$ | $2w_5-3$ | $2$ | $-2w_5$ |
| $f_{35}$ | $w_{17}-1$ | $-w_{17}$ | $1$ | $-1$ | $w_{17}$ | $w_{17}+2$ |
| $f_{39}$ | $\sqrt{2}-1$ | $1$ | $-2\sqrt{2}$ | $2\sqrt{2}$ | $-2$ | $-1$ |
| $f_{40}$ | $0$ | $2i$ | $-2i-1$ | $-2i$ | $-4$ | $4i$ |
| $f_{45}^{\mathrm{CM}}$ | $\sqrt{-5}$ | $0$ | $-\sqrt{-5}$ | $0$ | $0$ | $0$ |
| $f_{48}^{\mathrm{CM}}$ | $0$ | $2w_{-3}-1$ | $0$ | $-4w_{-3}+2$ | $0$ | $-2$ |
| $f_{52}$ | $0$ | $-w_{-3}+1$ | $0$ | $w_{-3}-2$ | $-3w_{-3}-3$ | $4w_{-3}-3$ |
| $f_{63}$ | $\sqrt{3}$ | $0$ | $-2\sqrt{3}$ | $1$ | $2\sqrt{3}$ | $2$ |
| $f_{64,A}^{\mathrm{CM}}$ | $0$ | $2i$ | $0$ | $0$ | $-6i$ | $0$ |
| $f_{64,B}$ | $0$ | $-i+1$ | $-i-1$ | $2i$ | $-i-1$ | $i-1$ |
| $f_{67}$ | $w_5-2$ | $-w_5-1$ | $-3$ | $3w_5-2$ | $-2w_5+1$ | $-3w_5-2$ |
| $f_{68}$ | $0$ | $\sqrt{3}+1$ | $-2\sqrt{3}$ | $-\sqrt{3}-1$ | $\sqrt{3}-3$ | $2\sqrt{3}+2$ |
| $f_{73}$ | $w_5-2$ | $-w_5-1$ | $w_5-2$ | $-3$ | $-w_5-1$ | $3w_5-1$ |
| $f_{80}$ | $0$ | $2i$ | $2i-1$ | $-2i$ | $4$ | $-4i$ |
| $f_{81}$ | $\sqrt{3}$ | $0$ | $-\sqrt{3}$ | $2$ | $-2\sqrt{3}$ | $-1$ |
| $f_{85}$ | $\sqrt{2}-1$ | $-\sqrt{2}-2$ | $-1$ | $\sqrt{2}-2$ | $\sqrt{2}-4$ | $-2\sqrt{2}$ |
| $f_{87}$ | $w_5$ | $1$ | $-2w_5+2$ | $-2w_5-1$ | $2w_5+1$ | $4w_5-3$ |
| $f_{88}$ | $0$ | $w_{17}$ | $-w_{17}+2$ | $-2w_{17}$ | $-1$ | $2w_{17}-2$ |
| $f_{93}$ | $w_5-2$ | $-1$ | $-2w_5-1$ | $2w_5-3$ | $2w_5-4$ | $2w_5-2$ |
| $f_{100}$ | $0$ | $2i$ | $0$ | $2i$ | $0$ | $-2i$ |
| $f_{103}$ | $w_5-2$ | $-1$ | $-w_5-1$ | $-1$ | $w_5-2$ | $3w_5-3$ |
| $f_{104}$ | $0$ | $w_{17}$ | $-w_{17}+2$ | $-w_{17}$ | $-2w_{17}$ | $1$ |
| $f_{107}$ | $w_5-1$ | $-w_5-1$ | $-w_5-1$ | $2w_5-3$ | $2w_5+1$ | $-6$ |
| $f_{112,A}$ | $0$ | $-w_{-3}+1$ | $-w_{-3}+2$ | $-2w_{-3}-1$ | $w_{-3}+1$ | $0$ |
| $f_{112,B}$ | $0$ | $w_{-3}-1$ | $-w_{-3}+2$ | $2w_{-3}+1$ | $-w_{-3}-1$ | $0$ |
| $f_{112,C}$ | $0$ | $w_{-3}$ | $3w_{-3}-3$ | $2w_{-3}+1$ | $-3w_{-3}$ | $2$ |
| $f_{115}$ | $w_5-2$ | $-1$ | $-1$ | $-2w_5$ | $2w_5-2$ | $2w_5-5$ |
| $f_{117}$ | $\sqrt{3}$ | $0$ | $0$ | $2$ | $-2\sqrt{3}$ | $1$ |
| $f_{125}$ | $w_5-1$ | $-w_5-1$ | $0$ | $-3$ | $-3$ | $3w_5-3$ |
| $f_{128}^{\mathrm{CM}}$ | $0$ | $2\sqrt{-2}$ | $0$ | $0$ | $2\sqrt{-2}$ | $0$ |
| $f_{133}$ | $w_5-2$ | $w_5-2$ | $-2w_5+1$ | $-1$ | $w_5-5$ | $1$ |
| $f_{135}$ | $w_{13}$ | $0$ | $-1$ | $-2w_{13}+2$ | $-2w_{13}$ | $2w_{13}+2$ |
| $f_{136}$ | $0$ | $2w_5-2$ | $2$ | $-2w_5+2$ | $-2w_5+2$ | $4w_5-2$ |
| $f_{147}$ | $\sqrt{2}-1$ | $-1$ | $-\sqrt{2}-2$ | $0$ | $-2$ | $\sqrt{2}-4$ |
| $f_{148}$ | $0$ | $-w_{-3}$ | $-w_{-3}+2$ | $-w_{-3}$ | $0$ | $-w_{-3}+2$ |
| $f_{160,A}$ | $0$ | $2\sqrt{2}$ | $1$ | $-2\sqrt{2}$ | $-4\sqrt{2}$ | $-2$ |
| $f_{160,B}$ | $0$ | $i+1$ | $-2i+1$ | $-i+1$ | $6i$ | $i-1$ |
| $f_{160,C}$ | $0$ | $i+1$ | $2i+1$ | $3i-3$ | $-2i$ | $-3i+3$ |
| $f_{160,D}$ | $0$ | $-i-1$ | $2i+1$ | $-3i+3$ | $2i$ | $-3i+3$ |

| $f$ | $a_2$ | $a_3$ | $a_5$ | $a_7$ | $a_{11}$ | $a_{13}$ |
|---|---|---|---|---|---|---|
| $f_{160,E}$ | $0$ | $-i-1$ | $-2i+1$ | $i-1$ | $-6i$ | $i-1$ |
| $f_{161}$ | $w_5-1$ | $-1$ | $-2w_5$ | $-1$ | $4w_5-2$ | $2w_5-3$ |
| $f_{165}$ | $\sqrt{2}-1$ | $-1$ | $-1$ | $-2\sqrt{2}-2$ | $-1$ | $4\sqrt{2}$ |
| $f_{167}$ | $w_5-1$ | $-w_5$ | $-1$ | $w_5-3$ | $0$ | $-w_5-2$ |
| $f_{175}$ | $w_5$ | $2w_5-2$ | $0$ | $1$ | $2w_5+1$ | $-2w_5$ |
| $f_{176}$ | $0$ | $w_{17}-1$ | $w_{17}+1$ | $-2w_{17}+2$ | $1$ | $-2w_{17}$ |
| $f_{177}$ | $w_5-1$ | $-1$ | $-2w_5+1$ | $w_5-4$ | $2w_5-1$ | $-2w_5-3$ |
| $f_{184}$ | $0$ | $w_{17}-1$ | $2$ | $0$ | $-2w_{17}+2$ | $-w_{17}+3$ |
| $f_{188,A}$ | $0$ | $w_5-2$ | $-2w_5$ | $-w_5-3$ | $4w_5-4$ | $4w_5-4$ |
| $f_{188,B}$ | $0$ | $w_{13}$ | $0$ | $-w_{13}+3$ | $-2w_{13}+2$ | $2$ |
| $f_{189}$ | $\sqrt{3}$ | $0$ | $\sqrt{3}$ | $1$ | $-\sqrt{3}$ | $2$ |
| $f_{191}$ | $w_5-1$ | $-1$ | $-w_5$ | $-w_5$ | $w_5-1$ | $3w_5-5$ |
| $f_{192}^{\mathrm{CM}}$ | $0$ | $2w_{-3}-1$ | $0$ | $4w_{-3}-2$ | $0$ | $2$ |
| $f_{205}$ | $w_5-1$ | $-1$ | $-1$ | $-3w_5+3$ | $2w_5-5$ | $3w_5-3$ |
| $f_{207}$ | $\sqrt{2}-1$ | $0$ | $-\sqrt{2}-2$ | $\sqrt{2}-2$ | $-2\sqrt{2}$ | $0$ |
| $f_{208,A}$ | $0$ | $w_{-3}-1$ | $0$ | $-w_{-3}+2$ | $3w_{-3}+3$ | $4w_{-3}-3$ |
| $f_{208,B}$ | $0$ | $w_{17}-1$ | $w_{17}+1$ | $-w_{17}+1$ | $-2w_{17}+2$ | $1$ |
| $f_{209}$ | $\sqrt{2}$ | $-\sqrt{2}-1$ | $-1$ | $-\sqrt{2}-2$ | $-1$ | $3\sqrt{2}-2$ |
| $f_{213}$ | $w_5-1$ | $-1$ | $-w_5+1$ | $-3$ | $-2w_5-1$ | $3w_5-4$ |
| $f_{221}$ | $w_5-1$ | $w_5-2$ | $-2w_5+1$ | $-w_5$ | $3w_5-3$ | $-1$ |
| $f_{224,A}$ | $0$ | $2w_5$ | $-2w_5+2$ | $-1$ | $-4w_5+4$ | $2w_5+2$ |
| $f_{224,B}$ | $0$ | $2w_5-2$ | $2w_5$ | $1$ | $-4w_5$ | $-2w_5+4$ |
| $f_{243}$ | $\sqrt{3}$ | $0$ | $2\sqrt{3}$ | $-1$ | $-2\sqrt{3}$ | $5$ |
| $f_{256}^{\mathrm{CM}}$ | $0$ | $2\sqrt{2}$ | $0$ | $0$ | $-2\sqrt{2}$ | $0$ |
| $f_{261}$ | $w_5-1$ | $0$ | $-2$ | $-2w_5+1$ | $-2w_5-3$ | $4w_5-3$ |
| $f_{272,A}$ | $0$ | $2w_5$ | $2$ | $-2w_5$ | $-2w_5$ | $-4w_5+2$ |
| $f_{272,B}$ | $0$ | $\sqrt{3}-1$ | $2\sqrt{3}$ | $-\sqrt{3}+1$ | $\sqrt{3}+3$ | $-2\sqrt{3}+2$ |
| $f_{275}$ | $w_{13}-1$ | $-w_{13}$ | $0$ | $w_{13}-3$ | $-1$ | $-5$ |
| $f_{280}$ | $0$ | $w_{17}$ | $1$ | $1$ | $-w_{17}$ | $-3w_{17}+2$ |
| $f_{287}$ | $w_5-1$ | $-w_5$ | $w_5$ | $-1$ | $-1$ | $-2w_5-3$ |
| $f_{297}$ | $\sqrt{3}-1$ | $0$ | $-\sqrt{3}-1$ | $\sqrt{3}-2$ | $-1$ | $-\sqrt{3}-2$ |
| $f_{299}$ | $w_5$ | $-w_5$ | $-w_5-1$ | $-1$ | $w_5-2$ | $-1$ |
| $f_{315}$ | $\sqrt{2}-1$ | $0$ | $-1$ | $-1$ | $-2\sqrt{2}-2$ | $2\sqrt{2}-2$ |
| $f_{320}$ | $0$ | $2\sqrt{2}$ | $-1$ | $2\sqrt{2}$ | $-4\sqrt{2}$ | $2$ |
| $f_{351}$ | $w_5-1$ | $0$ | $w_5-2$ | $-2w_5+1$ | $-3w_5-1$ | $-1$ |
| $f_{357}$ | $\sqrt{2}$ | $-1$ | $-\sqrt{2}-1$ | $-1$ | $1$ | $-\sqrt{2}-3$ |
| $f_{368}$ | $0$ | $w_{17}$ | $2$ | $0$ | $-2w_{17}$ | $w_{17}+2$ |
| $f_{376,A}$ | $0$ | $w_5-1$ | $-2w_5$ | $w_5-2$ | $-2w_5-2$ | $2w_5$ |
| $f_{376,B}$ | $0$ | $w_5-1$ | $-2$ | $-3w_5+2$ | $-2$ | $-2w_5-2$ |
| $f_{380}$ | $0$ | $\sqrt{3}+1$ | $1$ | $2$ | $-2\sqrt{3}$ | $-\sqrt{3}-1$ |
| $f_{400}$ | $0$ | $2i$ | $0$ | $2i$ | $0$ | $2i$ |
| $f_{416,A}$ | $0$ | $w_{17}-1$ | $-w_{17}-1$ | $-w_{17}-1$ | $-2$ | $-1$ |
| $f_{416,B}$ | $0$ | $w_{17}$ | $w_{17}-2$ | $-w_{17}+2$ | $2$ | $-1$ |
| $f_{440,A}$ | $0$ | $w_{17}-1$ | $1$ | $w_{17}-1$ | $1$ | $2$ |
| $f_{440,B}$ | $0$ | $w_{17}$ | $-1$ | $-w_{17}+2$ | $1$ | $2w_{17}$ |
| $f_{448,A}$ | $0$ | $2w_5$ | $2w_5-2$ | $1$ | $-4w_5+4$ | $-2w_5-2$ |
| $f_{448,B}$ | $0$ | $2w_5-2$ | $-2w_5$ | $-1$ | $-4w_5$ | $2w_5-4$ |

| $f$ | $a_2$ | $a_3$ | $a_5$ | $a_7$ | $a_{11}$ | $a_{13}$ |
|---|---|---|---|---|---|---|
| $f_{476,A}$ | $0$ | $w_5-1$ | $-w_5$ | $-1$ | $-2w_5-2$ | $-2w_5$ |
| $f_{476,B}$ | $0$ | $w_{13}$ | $w_{13}-1$ | $1$ | $0$ | $-2w_{13}+4$ |
| $f_{512,A}$ | $0$ | $\sqrt{2}$ | $-2$ | $-2\sqrt{2}$ | $-3\sqrt{2}$ | $-6$ |
| $f_{512,B}$ | $0$ | $\sqrt{2}$ | $2$ | $2\sqrt{2}$ | $-3\sqrt{2}$ | $6$ |
| $f_{525}$ | $w_{13}-1$ | $-1$ | $0$ | $-1$ | $-3$ | $-2w_{13}$ |
| $f_{544}$ | $0$ | $\sqrt{2}$ | $-2$ | $-3\sqrt{2}$ | $\sqrt{2}$ | $-4$ |
| $f_{560}$ | $0$ | $w_{17}-1$ | $1$ | $-1$ | $-w_{17}+1$ | $3w_{17}-1$ |
| $f_{592}$ | $0$ | $w_{-3}$ | $-w_{-3}+2$ | $w_{-3}$ | $0$ | $-w_{-3}+2$ |
| $f_{621}$ | $w_5$ | $0$ | $-w_5$ | $-2w_5-1$ | $w_5-2$ | $3w_5-3$ |
| $f_{640,A}$ | $0$ | $2w_5$ | $-1$ | $-2w_5+2$ | $2$ | $4w_5-2$ |
| $f_{640,B}$ | $0$ | $2w_5-2$ | $-1$ | $-2w_5$ | $-2$ | $-4w_5+2$ |
| $f_{645}$ | $\sqrt{2}$ | $-1$ | $-1$ | $\sqrt{2}$ | $-\sqrt{2}-1$ | $-3\sqrt{2}-1$ |
| $f_{704,A}$ | $0$ | $w_{17}$ | $w_{17}-2$ | $0$ | $1$ | $-2$ |
| $f_{704,B}$ | $0$ | $w_{17}-1$ | $-w_{17}-1$ | $0$ | $-1$ | $-2$ |
| $f_{752,A}$ | $0$ | $w_{13}-1$ | $0$ | $-w_{13}-2$ | $-2w_{13}$ | $2$ |
| $f_{752,B}$ | $0$ | $w_5+1$ | $2w_5-2$ | $-w_5+4$ | $4w_5$ | $-4w_5$ |
| $f_{752,C}$ | $0$ | $w_5$ | $2w_5-2$ | $w_5+1$ | $-2w_5+4$ | $-2w_5+2$ |
| $f_{752,D}$ | $0$ | $w_5$ | $-2$ | $-3w_5+1$ | $2$ | $2w_5-4$ |
| $f_{768,A}$ | $0$ | $\sqrt{-2}-1$ | $-2\sqrt{-2}$ | $2\sqrt{-2}$ | $-2$ | $4$ |
| $f_{768,B}$ | $0$ | $\sqrt{-2}+1$ | $2\sqrt{-2}$ | $2\sqrt{-2}$ | $2$ | $4$ |
| $f_{783}$ | $w_5-1$ | $0$ | $1$ | $w_5-2$ | $-2w_5$ | $w_5-3$ |
| $f_{880,A}$ | $0$ | $w_{17}$ | $1$ | $w_{17}$ | $-1$ | $2$ |
| $f_{880,B}$ | $0$ | $w_{17}-1$ | $-1$ | $-w_{17}-1$ | $-1$ | $-2w_{17}+2$ |
| $f_{928}$ | $0$ | $i-1$ | $0$ | $-2i$ | $i-1$ | $-4i$ |
| $f_{928}$ | $0$ | $-i+1$ | $0$ | $2i$ | $-i+1$ | $-4i$ |
| $f_{952}$ | $0$ | $w_{13}$ | $-w_{13}-1$ | $-1$ | $-2w_{13}+2$ | $-4$ |
| $f_{1053}$ | $w_5-1$ | $0$ | $w_5+1$ | $w_5-2$ | $-3w_5-1$ | $-1$ |
| $f_{1088}$ | $0$ | $\sqrt{2}$ | $2$ | $3\sqrt{2}$ | $\sqrt{2}$ | $4$ |
| $f_{1120,A}$ | $0$ | $-w_{17}$ | $-1$ | $1$ | $w_{17}-4$ | $w_{17}-2$ |
| $f_{1120,B}$ | $0$ | $w_{17}$ | $-1$ | $-1$ | $-w_{17}+4$ | $w_{17}-2$ |
| $f_{1280,A}$ | $0$ | $\sqrt{2}$ | $-1$ | $-\sqrt{2}$ | $-2\sqrt{2}$ | $-2$ |
| $f_{1280,B}$ | $0$ | $\sqrt{2}$ | $1$ | $\sqrt{2}$ | $-2\sqrt{2}$ | $2$ |
| $f_{1280,C}$ | $0$ | $\sqrt{2}$ | $-1$ | $3\sqrt{2}$ | $2\sqrt{2}$ | $6$ |
| $f_{1280,D}$ | $0$ | $\sqrt{2}$ | $1$ | $-3\sqrt{2}$ | $2\sqrt{2}$ | $-6$ |
| $f_{1280,E}$ | $0$ | $\sqrt{2}$ | $1$ | $3\sqrt{2}$ | $4\sqrt{2}$ | $-2$ |
| $f_{1280,F}$ | $0$ | $\sqrt{2}$ | $-1$ | $-3\sqrt{2}$ | $4\sqrt{2}$ | $2$ |
| $f_{1280,G}$ | $0$ | $\sqrt{3}-1$ | $-1$ | $\sqrt{3}+1$ | $-2$ | $-2\sqrt{3}$ |
| $f_{1280,H}$ | $0$ | $\sqrt{3}+1$ | $-1$ | $\sqrt{3}-1$ | $2$ | $2\sqrt{3}$ |
| $f_{1312}$ | $0$ | $\sqrt{2}$ | $0$ | $-\sqrt{2}$ | $-\sqrt{2}$ | $-2$ |
| $f_{1520}$ | $0$ | $-\sqrt{3}-1$ | $1$ | $-2$ | $2\sqrt{3}$ | $-\sqrt{3}-1$ |
| $f_{1792,A}$ | $0$ | $-2w_5+2$ | $2w_5-2$ | $1$ | $4$ | $-2w_5+2$ |
| $f_{1792,B}$ | $0$ | $2w_5-2$ | $2w_5-2$ | $-1$ | $-4$ | $-2w_5+2$ |
| $f_{1820}$ | $0$ | $w_{13}-1$ | $-1$ | $1$ | $-2$ | $-1$ |
| $f_{1904,A}$ | $0$ | $w_5$ | $w_5-1$ | $1$ | $-2w_5+4$ | $2w_5-2$ |
| $f_{1904,B}$ | $0$ | $-w_{13}$ | $w_{13}-1$ | $-1$ | $0$ | $-2w_{13}+4$ |
| $f_{1904,C}$ | $0$ | $w_{13}-1$ | $w_{13}-2$ | $1$ | $-2w_{13}$ | $-4$ |
| $f_{1916}$ | $0$ | $w_5-1$ | $w_5-1$ | $-w_5+1$ | $-w_5+1$ | $-2w_5$ |

| $f$ | $a_2$ | $a_3$ | $a_5$ | $a_7$ | $a_{11}$ | $a_{13}$ |
|---|---|---|---|---|---|---|
| $f_{2080}$ | 0 | $\sqrt{2}$ | $-1$ | 0 | $-3\sqrt{2}$ | 1 |
| $f_{2240,A}$ | 0 | $w_{17}$ | $-1$ | 1 | $w_{17}$ | $w_{17}-2$ |
| $f_{2240,B}$ | 0 | $-w_{17}$ | $-1$ | $-1$ | $-w_{17}$ | $w_{17}-2$ |
| $f_{2624}$ | 0 | $\sqrt{2}$ | 0 | $\sqrt{2}$ | $-\sqrt{2}$ | 2 |
| $f_{3159}$ | $w_5$ | 0 | $w_5-2$ | $2w_5-1$ | $-2$ | $-1$ |
| $f_{4160}$ | 0 | $\sqrt{2}$ | 1 | 0 | $-3\sqrt{2}$ | $-1$ |
| $f_{7280}$ | 0 | $w_{13}$ | $-1$ | $-1$ | 2 | $-1$ |
| $f_{7424,A}$ | 0 | $\sqrt{2}$ | 0 | $2\sqrt{2}$ | $\sqrt{2}$ | 2 |
| $f_{7424,B}$ | 0 | $\sqrt{2}$ | 0 | $-2\sqrt{2}$ | $\sqrt{2}$ | $-2$ |
| $f_{7664}$ | 0 | $-w_5+1$ | $w_5-1$ | $w_5-1$ | $w_5-1$ | $-2w_5$ |

## References

[AL78]   A. O. L. Atkin and W. C. W. Li, *Twists of newforms and pseudo-eigenvalues of W-operators*, Invent. Math. **48** (1978), no. 3, 221–243.

[BCDT99] C. Breuil, B. Conrad, F. Diamond, and R. Taylor, *On the modularity of elliptic curves over $\mathbb{Q}$: wild 3-adic exercises*, preprint (1999).

[BG97]   P. Bayer and J. González, *On the Hasse-Witt invariants of modular curves*, Experiment. Math. **6** (1997), no. 1, 57–76.

[Bol88]  O. Bolza, *On binary sextics with linear transformations into themselves*, Amer. J. Math. (1888), no. 10, 47–70.

[Bru95]  A. Brumer, *The rank of $J_0(N)$*, Astérisque (1995), no. 228, 3, 41–68, Columbia University Number Theory Seminar (New York, 1992).

[Car86]  H. Carayol, *Sur les représentations l-adiques associées aux formes modulaires de Hilbert*, Ann. Sci. École Norm. Sup. (4) **19** (1986), no. 3, 409–468.

[CGLR99] G. Cardona, J. González, J. C. Lario, and A. Rio, *On curves of genus 2 with Jacobian of $GL_2$-type*, Manuscripta Math. **98** (1999), no. 1, 37–54.

[Cre97]  J. E. Cremona, *Algorithms for modular elliptic curves*, second ed., Cambridge University Press, Cambridge, 1997.

[DS74]   P. Deligne and J-P. Serre, *Formes modulaires de poids* 1, Ann. Sci. École Norm. Sup. (4) **7** (1974), 507–530 (1975).

[FH99]   M. Furumoto and Y. Hasegawa, *Hyperelliptic quotients of modular curves $X_0(N)$*, Tokyo J. Math. **22** (1999), no. 1, 105–125.

[GL01]   J. González and J. C. Lario, *$\mathbb{Q}$-curves and their Manin ideals*, to appear in Amer. J. Math. (2001).

[Gon91]  J. González, *Equations of hyperelliptic modular curves*, Ann. Inst. Fourier (Grenoble) **41** (1991), no. 4, 779–795.

[Has97]  Y. Hasegawa, *Hyperelliptic modular curves $X_0^*(N)$*, Acta Arith. **81** (1997), no. 4, 369–385.

[HH96]   Y. Hasegawa and K. Hashimoto, *Hyperelliptic modular curves $X_0^*(N)$ with square-free levels*, Acta Arith. **77** (1996), no. 2, 179–193.

[IM91]   N. Ishii and F. Momose, *Hyperelliptic modular curves*, Tsukuba J. Math. **15** (1991), no. 2, 413–423.

[Liu94]  Q. Liu, *Conducteur et discriminant minimal de courbes de genre* 2, Compositio Math. **94** (1994), no. 1, 51–79.

[Mes81]  J.-F. Mestre, *Corps euclidiens, unités exceptionnelles et courbes élliptiques*, J. Number Theory **13** (1981), no. 2, 123–137.

[Mil72]  J. S. Milne, *On the arithmetic of abelian varieties*, Invent. Math. **17** (1972), 177–190.

[Mom81]  F. Momose, *On the l-adic representations attached to modular forms*, J. Fac. Sci. Univ. Tokyo Sect. IA Math. **28** (1981), no. 1, 89–109.

[Mur92]  N. Murabayashi, *On normal forms of modular curves of genus* 2, Osaka J. Math. **29** (1992), no. 2, 405–418.

[Ogg74]  A. P. Ogg, *Hyperelliptic modular curves*, Bull. Soc. Math. France **102** (1974), 449–462.

[Rib77]  K. A. Ribet, *Galois representations attached to eigenforms with Nebentypus*, 17–51. Lecture Notes in Math., Vol. 601.
[Rib80]  K. A. Ribet, *Twists of modular forms and endomorphisms of abelian varieties*, Math. Ann. **253** (1980), no. 1, 43–62.
[Shi71]  G. Shimura, *Introduction to the arithmetic theory of automorphic functions*, Publications of the Mathematical Society of Japan, No. 11. Iwanami Shoten, Publishers, Tokyo, 1971, Kanô Memorial Lectures, No. 1.
[Shi95]  M. Shimura, *Defining equations of modular curves $X_0(N)$*, Tokyo J. Math. **18** (1995), no. 2, 443–456.
[ST61]   G. Shimura and Y. Taniyama, *Complex multiplication of abelian varieties and its applications to number theory*, The Mathematical Society of Japan, Tokyo, 1961.
[Ste99]  W. A. Stein, *Hecke: The modular forms calculator*, Software (available online) (1999).

Department de Matemàtiques, Universitat Autònoma de Barcelona, Bellaterra, Barcelona, E-08193, Spain
 *E-mail address*: `enrikegj@mat.uab.es`

Escola Universitària Politècnica de Vilanova i la Geltrú, Av. Víctor Balaguer s/n, E-08800 Vilanova i la Geltrú, Spain
 *E-mail address*: `josepg@mat.upc.es`